\documentclass[12pt]{amsart}
\usepackage{amsfonts,amsmath, amsthm, amscd, amssymb}
\usepackage[hypertex]{hyperref}

\newtheorem{thm}{Theorem}
\newtheorem{lem}[thm]{Lemma}
\newtheorem{corollary}{Corollary}
\newtheorem{rmk}{Remark}
\newtheorem{example}{Example}
\newtheorem{defn}{Definition}

\newcommand{\cal}[1]{\mathcal #1}

\begin{document}

\title{Membership Problem in groups acting freely on $\mathbb{Z}^n$-trees}

\author{\textsf{Andrey Nikolaev} \and \textsf{Denis Serbin}}
\date{}
\address{Stevens Institute of Technology, Department of Mathematical Sciences,\\
1 Castle Point on Hudson, Hoboken, NJ 07030}

\begin{abstract}
Groups acting freely on $\mathbb Z^n$-trees ($\mathbb Z^n$-free groups) play a key role in the study
of non-archimedean group actions. Following Stallings' ideas, we develop graph-theoretic techniques
to investigate subgroup structure of $\mathbb Z^n$-free groups. As an immediate application of the
presented method, we give an effective solution to the Uniform Membership Problem and the Power
Problem in $\mathbb Z^n$-free groups.
\end{abstract}


\maketitle


\tableofcontents

\section{Introduction}
\label{sec:1}
The goal of the present work is to develop combinatorial techniques for studying algorithmic problems
in groups acting freely on $\mathbb{Z}^n$-trees, most notably The Uniform Membership Problem.

\subsection{The Membership Problems in groups}
\label{sub:1.1}
Recall that a finitely generated group
$$G = \langle x_1, x_2, \ldots, x_k \mid r_1, r_2, \ldots \rangle$$
is said to have {\em solvable membership problem} (or {\em solvable uniform membership problem}) if
there is an algorithm which, for any finite family of words $u, w_1, w_2, \ldots, w_n$ in $\{x_1,
x_2, \ldots, x_k\}^{\pm 1}$ decides whether or not the element of $G$ represented by $u$ belongs to
the subgroup of $G$ generated by the elements of $G$ corresponding to $w_1, w_2, \ldots, w_n$ (this
definition does not depend on the choice of a finite generating set for $G$). Similarly, if $H
\leqslant G$ is a specific subgroup, then $H$ is said to have solvable membership problem in $G$ if
there is an algorithm deciding for any word $u$ in $\{x_1, x_2, \ldots, x_k\}^{\pm 1}$ whether $u$
represents an element of $H$.

The membership problem in groups (and subgroup structure in general) has been extensively studied for
many classes of groups, employing multitude of different techniques and ideas. It is outside the
scope of this paper to give a comprehensive account of results in this area, but we would like to
point out certain developments that are of particular interest to us.

The framework for our study comes from the seminal paper \cite{St1} by J. Stallings, where he
introduced an extremely useful notion of a folding of graphs and initiated the study of subgroups
(and automorphisms) of free groups via folded directed labeled graphs. This approach turned out to
be very influential, it allowed researchers to prove many new results and simplify old proofs.
Detailed exposition of these results can be found, for example, in \cite{KM}.

Later this method was used by A. Myasnikov, V. Remeslennikov, and D. Serbin who showed that elements
of $F^{\mathbb Z[t]}$ (and, hence, of all its subgroups which are fully residually free groups) can
be viewed as reduced infinite words over the generating set of $F$ (see \cite{MRS}). Then it turned
out that many algorithmic problems for finitely generated fully residually free groups can be solved
by essentially the same methods as in standard free groups. Indeed, in \cite{MRS2} an analog of
Stallings' foldings was introduced for an arbitrary finitely generated subgroup of $F^{\mathbb Z[t]}$,
which allows one to solve effectively the membership problem in $F^{\mathbb Z[t]}$, as well as in
its arbitrary finitely generated subgroup. Next, in \cite{KMRS2} this technique was further developed
to obtain a solution to many algorithmic problems, such as the intersection problem, the conjugacy
problem, the malnormality problem etc. In \cite{NS} this technique was used to obtain an algorithmic
solution to the finite index problem and several other results on subgroup structure of fully
residually free groups.

Groups acting freely on $\mathbb Z^n$-trees ($\mathbb Z^n$-free groups) generalize fully residually
free groups, and in the present work we attempt to initiate the investigation of their subgroup
structure employing graph-theoretic techniques similar to those mentioned above.

Note that $\mathbb Z^n$-free groups are hyperbolic relative to maximal abelian subgroups. The
subgroup membership problem is in general undecidable in hyperbolic groups, as shown by Rips \cite{R}.
One notable point is that in our study we present a large subclass of relatively hyperbolic groups
for which the membership problem is decidable.

The authors are extremely grateful to Alexei Myasnikov, Olga Kharlampovich and Volker Diekert for
insightful discussions and many helpful comments and suggestions.

\subsection{$\mathbb{Z}^n$-free groups}
\label{sub:1.2}

The class of groups acting freely on $\mathbb{Z}^n$-trees became very important in the study of
actions on $\Lambda$-trees, where $\Lambda$ is an arbitrary ordered abelian group. It turned out
that group actions on $\mathbb{Z}^n$-trees play a key role in building a unified theory of
non-Archimedean actions, length functions and infinite words. The natural effectiveness of all
constructions (which is not the case for $\mathbb R$-trees) came along with a robust algorithmic
theory.

For detailed exhibition of the context of $\mathbb Z^n$-free groups we refer the reader to \cite{KMRS}.
Here we limit ourselves to noting that groups with regular $\mathbb Z^n$-free action play a key role
in recent advances (through the effort of Chiswell and M\"{u}ller \cite{ChMu}, and Miasnikov,
Kharlampovich, Remeslennikov and Serbin \cite{KMRS, KMS}) towards the solution on the following
long-standing problem of Alperin--Bass in the non-archimedean case.

\medskip

{\bf Problem.} {\it Describe finitely presented (finitely generated) groups acting freely on an
arbitrary $\Lambda$-tree.}

\medskip

Characterization of finitely generated groups with free regular length functions in $\mathbb{Z}^n$
given in \cite{KMRS} is a departure point of our present work, and we proceed by introducing the
necessary terminology, citing relevant results and inspecting class of finitely generated groups with
free regular length functions in $\mathbb{Z}^n$ in more detail.

\section{Preliminaries}
\label{sec:1+}

The notions introduced below can be found in a number of other works, for example, in \cite{KMS,
KMRS}. However, to keep the present exposition reasonably self-contained we choose to include a
number of definitions and statements here.

\subsection{Lyndon length functions and free actions on trees}
\label{sub:Lyndon}

Let $G$ be a group and $A$ an ordered abelian group. Then a function $l: G \rightarrow A$ is called
a {\it (Lyndon) length function} on $G$ if the following conditions hold:
\begin{enumerate}
\item[(L1)] $\forall\ g \in G:\ l(g) \geqslant 0$ and $l(1) = 0$,
\item[(L2)] $\forall\ g \in G:\ l(g) = l(g^{-1})$,
\item[(L3)] $\forall\ g, f, h \in G:\ c(g,f) > c(g,h)
\rightarrow c(g,h) = c(f,h)$,

\noindent where $c(g,f) = \frac{1}{2}(l(g)+l(f)-l(g^{-1}f))$.
\end{enumerate}

Notice that  $c(g,f)$ may not be defined in $A$ (if $l(g) + l(f) - l(g^{-1}f)$ is not divisible by
$2$), so in the axiom (L3) we assume that $A$ is canonically embedded into a divisible ordered
abelian group $A_{\mathbb{Q}} = A \otimes_\mathbb{Z} \mathbb{Q}$ (see \cite{MRS} for details).

It is not difficult to derive the following two properties of length functions from the axioms
(L1)--(L3):
\begin{itemize}
\item $\forall\ g, f \in G:\ l(g f) \leqslant l(g) + l(f)$,
\item $\forall\ g, f \in G:\ 0 \leqslant c(g,f) \leqslant \min\{l(g),l(f)\}$.
\end{itemize}
A length function $l:G \rightarrow A$ is called {\em free} if it satisfies the following two axioms.
\begin{enumerate}
\item[(L4)] $\forall\ g \in G:\ c(g,f) \in A$,
\item[(L5)] $\forall\ g \in G:\ g \neq 1 \rightarrow l(g^2) > l(g)$.
\end{enumerate}

For elements $g_1, \ldots, g_n \in G$ we write
$$g = g_1 \circ \cdots \circ g_n$$
if $g = g_1 \cdots g_n$ and $l(g) = l(g_1) + \cdots + l(g_n)$.

An {\em  $A$-metric space} is a pair $(X, d)$, where $d: X \times X \to A$ is a function such that
\begin{enumerate}
\item[(a)] $d(x, y) = 0 \iff x = y$,
\item[(b)] $d(x, y) \leqslant d(x, z) + d(y, z)$ for all $x, y, z \in X$.
\end{enumerate}
In this event $d$ is called {\em distance function}. A basic example of an $A$-metric space is $A$
itself with the usual distance function. An {\em isometry} between $A$-metric spaces is a map that
preserves distance.

Recall that for elements $a, b \in A$ the {\it closed segment} $[a,b]$ is defined as
$$[a,b] = \{c \in A \mid a \leqslant c \leqslant b \}.$$
Let $X$ be an $A$-metric space. A {\em closed segment} in $X$ is an isometry $\gamma : [0, a] \to X$.
If for every points $x,y \in X$ there exists a unique segment $\gamma$ such that $\gamma(0)$ and
$\gamma(a)$ then $X$ is said to be {\em geodesically convex}. The image of such $\gamma$ is denoted
$[x,y]$.

The following definition is due to Morgan and Shalen \cite{MoSh}.
\begin{defn}
An {\em $A$-tree} is a geodesically convex $A$-metric space $(X,d)$ such that for all $x,y,z \in X$
\begin{itemize}
\item $[x,y] \cap [y,z] = [y,w]$ for some $w \in X$,
\item $[x,y] \cap [y,z] = \{ y \} \Rightarrow [x,z] = [x,y] \cup [y,z]$.
\end{itemize}
\end{defn}
In the case $A = \mathbb{Z}$ this definition coincides with the definition of a simplicial tree with
the path metric. Note that $A$ viewed as an $A$-metric space is an $A$-tree.

\begin{defn}
An isometric action of a group on an $A$-tree $X$ is free if there are no inversions and the
stabilizer of each point of $X$ is trivial. We say that a group $G$ is $A$-free if $G$ admits such an
action on an $A$-tree.
\end{defn}

The following theorem due to Chiswell provides a link between Lyndon length function and actions on
$A$-trees.
\begin{thm}
\label{chis0}\cite{Ch1}
A group $G$ has a free Lyndon length function with values in $A$ if and only if $G$ acts freely on an
$A$-tree.
\end{thm}

\subsection{Regular length functions}
\label{sub:Lyndon-actions}

In this section we define regular length functions, and show some examples of groups with regular
length functions. The Theorem \ref{main4} gives a host of new examples of such groups.

A length function $l: G \rightarrow A$ is called {\it regular} if it satisfies the {\it regularity}
axiom:
\begin{enumerate}
\item[(L6)] $\forall\ g, f \in G,\ \exists\ u, g_1, f_1 \in G:$
$$g = u \circ g_1 \ \& \  f = u \circ f_1 \ \& \ l(u) = c(g,f).$$
\end{enumerate}

Here are several examples of groups with regular free length functions.

\begin{example}
\label{ex:1}
Let $F = F(X)$ be a free group on $X$. The length function
$$|\ | : F \rightarrow \mathbb{Z},$$
where $|f|$ is a natural length of $f \in F$ as a finite word, is regular since the common initial
subword of any two elements of $F$ always exists and belongs to $F$.
\end{example}

\begin{example}
\label{ex:2}
In \cite{MRS} it was proved that the Lyndon's free $\mathbb{Z}[t]$-group has a regular free length
function with values in $\mathbb{Z}[t]$.
\end{example}

\begin{example}
\label{ex:3} \cite{KMRS}
Let $F = F(X)$ be a free group on $X$. Consider an HNN-extension
$$G = \langle F, s \mid u^s = v \rangle,$$
where $u,v \in F$ are such that $|u| = |v|$ and $u$ is not conjugate to $v^{-1}$. Then there is a
regular free length function $l : G \rightarrow \mathbb{Z}^2$ which extends the natural
integer-valued length function on $F$.
\end{example}

For more involved examples, we refer the reader to \cite{KMRS}.

\subsection{Infinite words and length functions}
\label{subs:words}

In this subsection at first we recall some notions from the theory of ordered abelian groups (for
all the details we refer the reader to the books \cite{Gl} and \cite{KopMed}) and then following
\cite{MRS} describe the construction of infinite words.

\smallskip

An ordered abelian group $A$ is called {\it discretely ordered} if $A$ has a minimal positive element
(we denote it by $1_A$). In this event, for any $a \in A$ the following hold:
\begin{enumerate}
\item[(1)] $a + 1_A = \min\{b \mid b > a\}$,
\item[(2)] $a - 1_A = \max\{b \mid b < a\}$.
\end{enumerate}

Observe that if $A$ is any ordered abelian group then $\mathbb{Z} \oplus A$ is discretely ordered
with respect to the right lexicographic order.

\smallskip

Let $A$ be a discretely ordered abelian group and let $X = \{x_i \mid i \in I\}$ be a set. Put
$X^{-1} = \{x_i^{-1} \mid i \in I\}$ and $X^\pm = X \cup X^{-1}$. An {\em $A$-word} is a function of
the type
$$w: [1_A, \alpha_w] \to X^\pm,$$
where $\alpha_w \in A,\ \alpha_w \geqslant 0$. The element $\alpha_w$ is called the {\em length}
$|w|$ of $w$.

In the case $A=\mathbb Z^n$, represent $Z^n$ as a union of the chain of subgroups
$$\{0\} = A_0 < A_1 < A_2 < \cdots < A_n = \mathbb Z^n,$$
where
$$A_i = \{(a_1, a_2, \ldots, a_n) \mid a_{i+1} = a_{i+2} = \cdots = a_n = 0\} \simeq \mathbb{Z}^i.$$
Then we say that a word $w$ has {\em height} $i$ if $|w| \in A_i - A_{i-1}$, denoted $ht(w) = i$.
Note that in the case of an arbitrary discretely ordered $A$, an analogous definition can be given
in terms of convex subgroups of $A$. Height of an element of a group with a free $Z^n$-valued Lyndon
length function is defined in a similar fashion.

\smallskip

By $W(A,X)$ we denote the set of all $A$-words. Observe, that $W(A,X)$ contains an empty $A$-word
which we denote by $\varepsilon$.

Concatenation $u v$ of two $A$-words $u,v \in W(A,X)$ is an $A$-word of length $|u| + |v|$ and such
that:
\[(uv)(a) = \left\{ \begin{array}{ll}
\mbox{$u(a)$} & \mbox{if $1_A \leqslant a \leqslant |u|$} \\
\mbox{$v(a - |u|)$ } & \mbox{if $|u|  < a \leqslant |u| + |v|$}
\end{array}
\right. \]

An $A$-word $w$ is {\it reduced} if $w(\beta + 1_A) \neq w(\beta)^{-1}$ for each $1_A \leqslant
\beta < |w|$.  If the concatenation $uv$ of two reduced $A$-words $u$ and $v$ is also reduced then
we write $u v = u \circ v$. We denote by $R(A,X)$ the set of all reduced $A$-words. Clearly,
$\varepsilon \in R(A,X)$.

\smallskip

For $u \in W(A,X)$ and $\beta \in [1_A, \alpha_u]$ by $u_\beta$ we denote the restriction of $u$ on
$[1_A,\beta]$. If $u \in R(A,X)$ and $\beta \in [1_A, \alpha_u]$ then
$$u = u_\beta \circ {\tilde u}_\beta,$$
for some uniquely defined ${\tilde u}_\beta$.

An element ${\rm com}(u,v) \in R(A,X)$ is called the ({\emph{longest}) {\it common initial segment}
of $A$-words $u$ and $v$ if
$$u = {\rm com}(u,v) \circ \tilde{u}, \ \ v = {\rm com}(u,v) \circ \tilde{v}$$
for some (uniquely defined) $A$-words $\tilde{u}, \tilde{v}$ such that $\tilde{u}(1_A) \neq
\tilde{v}(1_A)$.

Now, we can define the product of two $A$-words. Let $u,v \in R(A,X)$. If ${\rm com}(u^{-1}, v)$ is
defined then
$$u^{-1} = {\rm com}(u^{-1},v) \circ {\tilde u}, \ \ v = {\rm com} (u^{-1},v) \circ {\tilde v},$$
for some uniquely defined ${\tilde u}$ and ${\tilde v}$. In this event put
$$u \ast v = {\tilde u}^{-1} \circ {\tilde v}.$$
The product ${\ast}$ is a partial binary operation on $R(A,X)$.

\smallskip

An element $v \in R(A,X)$ is termed {\it cyclically reduced} if $v(1_A)^{-1} \neq v(|v|)$. We say
that an element $v \in R(A,X)$ admits a {\it cyclic decomposition} if $v = c^{-1} \circ u \circ c$,
where $c, u \in R(A,X)$ and $u$ is cyclically reduced. Observe that a cyclic decomposition is unique
(whenever it exists). We denote by $CR(A,X)$ the set of all cyclically reduced words in $R(A,X)$ and
by $CDR(A,X)$ the set of all words from $R(A,X)$ which admit a cyclic decomposition.

\smallskip

Below we refer to $A$-words as {\it infinite words} usually omitting $A$ whenever it does not produce
any ambiguity.

The following result establishes the connection between infinite words and length functions.
\begin{thm}
\label{co:3.1} \cite{MRS}
Let $A$ be a discretely ordered abelian group and $X$ be a set. Then any subgroup $G$ of $CDR(A,X)$
has a free Lyndon length function with values in $A$ -- the restriction $L \mid_G$ on $G$ of the
standard length function $L$ on $CDR(A,X)$.
\end{thm}

The converse of the theorem above was obtained by I.Chiswell \cite{Ch}.

\begin{thm}
\label{chis} \cite{Ch}
Let $G$ have a free Lyndon length function $L : G \rightarrow A$, where $A$ is a discretely ordered
abelian group. Then there exists a length preserving embedding $\phi : G \rightarrow CDR(A,X)$, that
is, $|\phi(g)| = L(g)$ for any $g \in G$.
\end{thm}

\begin{corollary}
\label{chis-cor} \cite{Ch}
Let $G$ have a free Lyndon length function $L : G \rightarrow A$, where $A$ is an arbitrary ordered
abelian group. Then there exists an embedding $\phi : G \rightarrow CDR(A',X)$, where $A' = \mathbb{Z}
\oplus A$ is discretely ordered with respect to the right lexicographic order and $X$ is some set,
such that, $|\phi(g)| = (0,L(g))$ for any $g \in G$.
\end{corollary}

Theorems \ref{co:3.1} and \ref{chis}, and Corollary \ref{chis-cor} show that a group has a free Lyndon
length function if and only if it embeds into a set of infinite words and this embedding preserves
the length. Moreover, it is not hard to show that this embedding also preserves regularity of the
length function.

\begin{thm}
\label{chis-cor-1} \cite{KMS}
Let $G$ have a free regular Lyndon length function $L : G \rightarrow A$, where $A$ is an arbitrary
ordered abelian group. Then there exists an embedding $\phi : G \rightarrow R(A',X)$, where $A'$ is
a discretely ordered abelian group and $X$ is some set, such that, the Lyndon length function on
$\phi(G)$ induced from $R(A',X)$ is regular.
\end{thm}

\subsection{Groups with free length functions in $\mathbb{Z}^n$}
\label{sub:structure}

Let $G$ be a finitely generated group acting freely on a $\mathbb{Z}^n$-tree for some $n \in
\mathbb{N}$. Hence, by Corollary \ref{chis-cor} \cite{Ch}, $G$ embeds into $CDR(X, \mathbb{Z}^n)$
for some alphabet $X$. Then according to Corollary 4 \cite{KMS}, there exists a finite alphabet $Y$
and an embedding $\phi : G \rightarrow G'$, where $G'$ is a subgroup of $CDR(\mathbb{Z}^n, Y)$ with
a regular length function, such that $|g|_G = |\phi(g)|_{G'}$ for every $g \in G$. Since $G$ can be
viewed as a subgroup of $G'$ then in order to solve algorithmic problems for subgroups of $G$ it is
enough to solve them for subgroups of $G'$. So, we can assume $G = G'$ and fix it for the rest of
the text.

On the other hand, the structure of groups with regular free length functions in $\mathbb{Z}^n$ is
known (see \cite{KMRS} for all the details).

\begin{thm} [Theorem 7 \cite{KMRS}]
\label{th:main2}
Let $G$ be a finitely generated group with a regular free Lyndon length function in $\mathbb{Z}^n$.
Then $G$ can be represented as a union of a finite series of groups
$$G_{(1)} < G_{(2)} < \cdots < G_{(n)} = G,$$
where $G_{(1)}$ is a free group of finite rank, $G_{(i)}$ has a regular free Lyndon length function
in $\mathbb{Z}^i$, and
\begin{equation}
\label{eq:0}
G_{(i+1)} = \langle G_{(i)}, s_{i,1},\ \ldots,\ s_{i,k_i} \mid s_{i,j}^{-1}\ C_{i,j}\ s_{i,j} =
\phi_{i,j}(C_{i,j}) \rangle,
\end{equation}
where for each $j \in [1,k_i],\ C_{i,j}$ and $\phi_{i,j}(C_{i,j})$ are cyclically reduced centralizers
of $G_{(i)}$, $\phi_{i,j}$ is an isomorphism, and the following conditions are satisfied:
\begin{enumerate}
\item[(1)] $C_{i,j} = \langle c^{(i,j)}_1, \ldots, c^{(i,j)}_{m_{i,j}} \rangle,\ \phi_{i,j}(C_{i,j})
= \langle d^{(i,j)}_1, \ldots, d^{(i,j)}_{m_{i,j}} \rangle$, where $\phi_{i,j}(c^{(i,j)}_k) =
d^{(i,j)}_k,\ k \in [1,m_{i,j}]$ and
$$ht(c^{(i,j)}_k) = ht(d^{(i,j)}_k) < ht(d^{(i,j)}_{k+1}) = ht(c^{(i,j)}_{k+1}),\ k \in [1,m_{i,j}-1],$$
\item[(2)] $|\phi_{i,j}(w)| = |w|$ for any $w \in C_{i,j}$,
\item[(3)] $w$ is not conjugate to $\phi_{i,j}(w)^{-1}$ in $G_{(i)}$ for any $w \in C_{i,j}$,
\item[(4)] if $A, B \in \{C_{i,1}, \phi_{i,1}(C_{i,1}), \ldots, C_{i,k_i}, \phi_{i,k_i}(C_{i,k_i})\}$
then either $A = B$, or $A$ and $B$ are not conjugate in $G_{(i)}$,
\item[(5)] for any $j_0 \in [1,k_i]$, the centralizer $C_{i,j_0}$ can appear in the list
$$\{C_{i,j} \mid j \in [1,k_i]\} \cup \{\phi_{i,j}(C_{i,j}) \mid j \in [1,k_i]\}$$
not more than twice.
\end{enumerate}
\end{thm}

\begin{rmk}
\label{rem:main_th}
In addition, from the proof of Theorem 7 \cite{KMRS} it follows that
\begin{enumerate}
\item[(6)] $ht(s_{i,j}) = ht(s_{i,k}) = i + 1$ for any $j, k \in [1,k_i]$,
\item[(7)] if $C_{i,j} = C_{i,k}$ then either $c(s_{i,j}, s_{i,k}) = 0$, or $s_{i,j} = s_{i,k}$.
\end{enumerate}
\end{rmk}

The converse to the theorem above also holds.

\begin{thm}
\label{main4} \cite{KMRS}
Let $H$ be a group with a regular free length function in $\mathbb{Z}[t]$. Let $A$ and $B$ be
centralizers in $H$ whose elements are cyclically reduced and such that there exists an isomorphism
$\phi : A \rightarrow B$ with the following properties
\begin{enumerate}
\item $a$ is not conjugate to $\phi_i(a)^{-1}$ in $H$ for any $a \in A$,
\item $|\phi(a)| = |a|$ for any $a \in A$.
\end{enumerate}
Then the group
\begin{equation}
\label{eq:G}
G = \langle H, z \mid z^{-1} A z = B \rangle,
\end{equation}
has a regular free length function in $\mathbb{Z}[t]$ which extends the length function on $H$.
\end{thm}

Hence, we get a lot of information using the fact that every $\mathbb Z^n$-free group nicely embeds
into a $\mathbb Z^n$-free group with the regular underlying length function and we are going to use
this fact a lot.

Below we list other notable properties of $\mathbb Z^n$-free groups. Items \ref{item1}, \ref{item2},
\ref{item5}, \ref{item6} can be found in Martino and Rourke's survey \cite{MarR}.
\begin{enumerate}
\item \label{item1}
$\mathbb Z^n$-free groups are commutation transitive, and any abelian subgroup of a $\mathbb{Z}^n$-free
group is free abelian of rank at most $n$. It follows by simple inspection of action of commuting
elements. This property holds for all $A$-free groups for an arbitrary abelian group $A$, as observed
by Bass in \cite{Ba}.

\item \label{item2}
$\mathbb Z^n$-free groups are coherent (that is, finitely generated $\mathbb Z^n$-free groups are
finitely presented). While this was, likely, first observed in \cite{MarR}, it follows from the
results of \cite{KMW}.

\item \label{item4}
$\mathbb Z^n$-free groups are hyperbolic relative to maximal abelian subgroups. Indeed, the base
change functor (see Theorem 2.3 in \cite{MarR}) allows to obtain a free $\mathbb R^n$-action and then
the statement follows by the result of Guirardel \cite{G}.

\item \label{item5}
A finitely generated $\mathbb Z^n$-free group all of whose maximal abelian subgroups are cyclic is
word hyperbolic (as are all its finitely generated subgroups).

\item \label{item6}
The Word Problem is decidable in any $\mathbb Z^n$-free group.

\item \label{item7}
The class of $\mathbb{Z}^n$-free groups is closed under amalgamated free products along maximal cyclic
subgroups ($n$ is not preserved)  \cite{MarR2}.
\end{enumerate}

\label{difficuly}

As we are concerned with solving the Membership Problem for $\mathbb Z^n$-free groups, Item \ref{item2}
is of particular interest, since one can try to employ Kapovich--Miasnikov--Weidmann graphs of groups
technique to solve the membership problem. Unfortunately, in the case of $\mathbb Z^n$-free groups
this technique does not provide an effective solution because the algorithm presented in \cite{KMW}
requires solving the Subgroup Intersection Problem in vertex groups.

\subsection{Groups with free regular length functions in $\mathbb{Z}^n$, a refinement}
\label{sub:2.1}

In the next subsections we introduce normal forms for elements of $G$. But in order to do that we need
a slight refinement of Theorem \ref{th:main2} (we keep the notation introduced in Subsection
\ref{sub:structure}).

First of all, for a subgroup $H$ of $G$ let the {\em height} of $H$ be defined as
$$ht(H) = \max\{ ht(h) \mid h \in H\}.$$
Obviously $ht(H)$ exists and does not exceed $n$. Next, denote
$$T^{(i)} = \{s_{i,j} \mid j \in [1,k_i]\},\ T = \bigcup_{i=1}^{n-1} T^{(i)},$$
$${\cal C}^{(i)} = \{C_{i,j} \mid j \in [1,k_i]\} \cup \{D_{i,j} \mid j \in [1,k_i]\},$$
$${\cal C} = \bigcup_{i=1}^{n-1} {\cal C}^{(i)}.$$
Consider $C \in {\cal C}$. We fix the generating set for $C$ which appears in Theorem \ref{th:main2}
and consider $c \in C$, the generator of minimal height. Observe that every $\varepsilon \neq u \in
C$ has either $c$ or $c^{-1}$ as an initial and terminal subword. Hence, for a fixed $c$ we can
introduce an {\em orientation} on $C$ with respect to $c$ as follows. We can represent $C -
\{\varepsilon\}$ as a union $C^+ \cup C^-$, where every $u \in C^+$ has $c$ as an initial and terminal
subword, and every $u \in C^-$ has $c^{-1}$ as an initial and terminal subword. Observe that
$(C^+)^{-1} = C^-$. Obviously, there are only two possible orientations for $C$, the one with respect
to $c$, and another one with respect to $c^{-1}$.

Now take another $D \in {\cal C}$ such that $w^{-1} \ast C \ast w = D$ for some $w \in G$. If
$C$ is oriented already then we can choose an orientation for $D$ so that $w^{-1} \ast C^+ \ast w
\leqslant D^+$. In this case we say that $C$ and $D$ have the same orientation.

Since $g^{-1} \ast f \ast g \neq f^{-1}$ for any $\varepsilon \neq f,g \in G$ it follows that we can
orient all elements of ${\cal C}$ so that any two conjugate centralizers have the same orientation.
Let us fix such an orientation on ${\cal C}$.

Following \cite{KMRS}, for $C \in {\cal C}$ we call $t \in T^{\pm 1}$ {\em attached to $C$} if
$ht(t) > ht(C)$ and $ht(t^{-1} \ast C \ast t) = ht(C)$. If $t$ is attached to $C$ then we call $t$
{\em left} if $c^{-1} \ast t = c^{-1} \circ t$, and {\it right} if $c \ast t = c \circ t$ for any
$c \in C^+$.

Now, assume that $C \in {\cal C}^{(k)}$ and let $c \in C^+$ be the generator of maximal height.
Observe that if there exists a right(left)-attached to $C$ with respect to $c$ element then by Lemma
15 \cite{KMRS} there exists $D \in \bigcup_{i=1}^k {\cal C}^{(i)}$ and $d \in D^+$ which is the generator
of maximal height, such that $c$ is conjugate to $d$ in $G_{(k)}$ and $D$ does not have
right(left)-attached with respect to $d$ elements. Next, by Lemma 16 \cite{KMRS} it follows that there
is no element $g \in G_{(k)}$ which has any positive power of $d$ as an initial(terminal) subword.

Starting from the upper level $n-1$ and using the above argument, in every pair of centralizers $C,
D \in {\cal C}^{(k)}$ such that $t^{-1} C t = D$ for some $t \in T^{(k)}, k \in [1,n-1]$ we can
substitute $C$ by $w_C^{-1} C' w_C$ and $D$ by $w_D^{-1} D' w_D$, where $w_C, w_D$ are products of some
elements from $\bigcup_{i=1}^k T^{(i)}$ and $C', D' \in \bigcup_{i=1}^k {\cal C}^{(i)}$ are such that
any element $g \in G_{(k)}$ has neither any positive power of $c$ as an initial subword,
nor any positive power of $d$ as a terminal subword for the generators of maximal height $c \in C'^+,\
d \in D'^+$. In other words, without loss of generality we can assume that $C' = C$ and $D' = D$ from
the beginning. Finally, observe that by Theorem \ref{th:main2}, $t$ (note that $t \notin G_{(k)}$) is
the only element of $T^{(k)}$ right-attached to $C$ with respect to $c$ and left-attached to $D$ with
respect to $d$.

Now, we split ${\cal C}^{(i)}$ into conjugacy classes ${\cal C}^{(i)}_j,\ j \in [1,K_i(G)]$ in
$G_{(i+1)}$. Observe that it is possible that $A \in {\cal C}^{(i)}_j$ is a subgroup of $A' \in
{\cal C}^{(i')}_{j'}$ for some $i \leqslant i'$. Also, by the above assumption, all elements of
${\cal C}^{(i)}_j$ for $i\in [1,n-1], j \in [1,K_i(G)]$ have the same orientation. Moreover, any two
$A, B \in {\cal C}^{(i)}_j$ have the same height, and $w^{-1} \ast A \ast w = B,\ w \neq \varepsilon$
implies $ht(w) > ht(A)$. Denote by $T^{(i)}_j$ the corresponding subset of $T^{(i)}$.

Denote $ht({\cal C}^{(i)}_j) = ht(A)$, where $A \in {\cal C}^{(i)}_j,\ j \in [1,K_i(G)]$. Let
${\cal C}_1, {\cal C}_2, \ldots,$ ${\cal C}_{K(G)}$ be a reordering of the set of conjugacy classes
$$\{{\cal C}^{(1)}_1, \ldots, {\cal C}^{(1)}_{K_1(G)}, {\cal C}^{(2)}_1,\ \ldots,
{\cal C}^{(2)}_{K_2(G)}, \ldots, {\cal C}^{(n-1)}_1, \ldots, {\cal C}^{(n-1)}_{K_{n-1}(G)}\}$$
such that ${\cal C}^{(i)}_j$ has lower number than ${\cal C}^{(i')}_{j'}$ whenever $i < i'$, while
the order relation on the conjugacy classes within the same ${\cal C}^{(i)}$ is arbitrary. Assume that
the corresponding sets $T^{(i)}_j$ are reordered accordingly. Observe that
$${\cal C} = \bigcup_{i=1}^{K(G)} {\cal C}_i,\ \ \ T = \bigcup_{i=1}^{K(G)} T_i.$$
Now, using the notation just introduced we are ready to reformulate Theorem \ref{th:main2} pointing
out the properties of the new presentation of $G$ necessary in the further considerations.

\begin{thm}
\label{th:1}
Let $T_i$ and ${\cal C}_i$, where $i \in [1, K(G)]$, be defined as above. Then $G$ can be represented
as a union of a finite series of groups
\begin{equation}
\label{eq:2}
G_1 < G_2 < \cdots < G_{K(G)} = G,
\end{equation}
where $G_1$ is a free group of finite rank, and
$$G_{i+1} = \langle G_i, T_i \mid t^{-1} C_t t \stackrel{\phi_t}{=} D_t,\ t \in T_i,\ C_t, D_t \in
{\cal C}_i \rangle,$$
where $\phi_t: C_t \rightarrow D_t$ is an isomorphism and the following conditions are satisfied:
\begin{enumerate}
\item[(1)] each $G_i$ is a group with a free regular length function,
\item[(2)] $C_t$ and $D_t$ are cyclically reduced centralizers for any $t \in T_i$,
\item[(3)] $C_t = \langle c^{(t)}_1, \ldots, c^{(t)}_{k_t} \rangle,\ D_t = \langle d^{(t)}_1, \ldots,
d^{(t)}_{k_t} \rangle$, where all $c^{(t)}_i, d^{(t)}_i$ are cyclically reduced,
$$ht(c^{(t)}_i) = ht(d^{(t)}_i) < ht(d^{(t)}_{i+1}) = ht(c^{(t)}_{i+1}),\ i \in [1,k_t-1],$$
and $\phi(c^{(t)}_i) = d^{(t)}_i$,
\item[(4)] $|\phi_t(c)| = |c|$ for any $c \in C_t$,
\item[(5)] for any $s \in T_i$, the centralizer $C_s$ can appear in the list
$$\{C_t \mid t \in T_i\} \cup \{D_t \mid t \in T_i\}$$
not more than twice,
\item[(6)] $ht(s) = ht(t)$ for any $s, t \in T_i$, and if $C_s = C_t$ then either $c(s, t) = 0$, or
$s = t$,
\item[(7)] any element $g \in G_i$ may have only a bounded positive power of $c^{(t)}_{k_t}$ as an
initial subword and only a bounded positive power of $d^{(t)}_{k_t}$ as a terminal subword.
\end{enumerate}
\end{thm}
{\it Proof.} (1) -- (5) follow from Theorem \ref{th:main2}, (6) follows from Remark \ref{rem:main_th}, and
(7) follows from the discussion above.

\hfill $\square$

\subsection{Normal forms in  $G_{i+1}$}
\label{sub:2.2}

Fix $i \in [1,K(G)]$. From now on we are going to use the following notation. Let $ST_i \subset
\langle T_i \rangle$ be the set of elements $w \in \langle T_i \rangle$ such that there exist $C, D
\in {\cal C}_i$ with the property $w^{-1} \ast C \ast w = D \in {\cal C}_i$. Obviously, $\varepsilon
\in ST_i$ and $(ST_i)^{-1} = ST_i$. For every $w \in ST_i$, if $w^{-1} \ast C \ast w = D$ for $C,D
\in {\cal C}_i$ then denote $C_w = C,\ D_w = D$, and denote by $u_w,\ v_w$ the generators of $C_w$
and $D_w$ of maximal height such that $w^{-1} \ast u_w \ast w = v_w,\ u_w \ast w = u_w \circ w,\ w
\ast v_w = w \circ v_w$. It is easy to see that $C_{w^{-1}} = D_w,\ D_{w^{-1}} = C_w$ and $u_{w^{-1}}
= v_w^{-1},\ v_{w^{-1}} = v_w^{-1}$. Also, let
$$U_i = \{ u_t, v_t \mid t \in T_i \}.$$

According to Theorem \ref{th:1}, $G_{i+1}$ is obtained from $G_i$ as a multiple HNN-extension, so
every $g \in G_{i+1} - G_i$ has the following representation
\begin{equation}
\label{eq:3}
g = g_1 \ast w_1 \ast g_2 \ast \cdots \ast w_k \ast g_{k+1},
\end{equation}
where $g_j \in G_i,\ j \in [1,k+1],\ w_j \in ST_i,\ j \in [1,k]$, and if $D_{w_j} = C_{w_{j+1}}$
then $g_j \notin D_{w_j}$ for any $j \in [1,k-1]$.

\begin{lem}
\label{le:2.2.1}
Let $g \in G_i$ for some $i \in [1, K(G)]$. If $c \in C,\ d \in D$ are generators of maximal height
for arbitrary $C, D \in {\cal C}_i$ then there exists $N(g) > 0$ such that for every $r \geqslant N(g)$
\begin{enumerate}
\item[(a)] $c^r \ast g = c^{r-N(g)} \circ (c^{N(g)} \ast g)$,
\item[(b)] $c^r \ast g \ast d^r = c^{r-N(g)} \circ (c^{N(g)} \ast g \ast d^{N(g)}) \circ d^{r-N(g)}$,
provided $g \notin C$ if $c = d^{-1}$.
\end{enumerate}
\end{lem}
{\it Proof.} Follows from Theorem \ref{th:1} (7).

\hfill $\square$

By Lemma \ref{le:2.2.1} there exist $N_j, M_j \in \mathbb{Z},\ j \in [1,k]$ such that $M_j, N_j \geqslant
N(g_j)$ and (\ref{eq:3}) can be transformed into
\begin{equation}
\label{eq:4}
\begin{array}{rl}
g =& (g_1 \ast u_{w_1}^{N_1}) \circ (u_{w_1}^{-N_1} \ast w_1 \ast v_{w_1}^{-M_1}) \circ (v_{w_1}^{M_1}
\ast g_2 \ast u_{w_2}^{N_2}) \circ\\
&\phantom{a}\\
&\circ \cdots \circ (u_{w_k}^{-N_k} \ast w_k \ast v_{w_k}^{-M_k}) \circ (v_{w_k}^{M_k} \ast g_{k+1}).
\end{array}
\end{equation}
Since $N_j, M_j$ only satisfy the condition $M_j,N_j > N(g_j)$, and $g_j$ may have elements of
$D_{w_{j-1}}$ and $C_{w_j}$ respectively as initial and terminal subwords it follows that (\ref{eq:4})
is not a unique representation of $g$.

Our goal is to obtain a unique representation for $g$.

Obviously $u^N \ast w \ast v^M = w \ast v^{M + N},\ (u^N \ast w \ast v^M)^{-1} = w^{-1} \ast u^{-(M
+ N)}$ for any $w \in ST_i,\ N, M \in \mathbb{Z}$, and $u = u_w^\delta,\ v = v_w^\delta,\ \delta \in
\{-1,1\}$, so, denote
$$w(r) = w \ast v^r.$$
Notice that $w(r) \in G_{i+1}$ for every $r \in \mathbb{Z}$, and $w(r)$ has any natural powers of
$u_w$ and $v_w$ respectively as an initial and terminal subword. So, we can use the notation
$C_{w(r)} = C_w,\ D_{w(r)} = D_w,\ u_{w(r)} = u_w,\ v_{w(r)} = v_w$.

Now, let
\begin{equation}
\label{eq:5}
g = g_1 \circ w_1(r_1) \circ g_2 \circ \cdots \circ w_k(r_k) \circ g_{k+1},
\end{equation}
where $g_j \in G_i,\ j \in [1,k+1],\ w_j \in ST_i,\ j \in [1,k]$ and if $D_{w_j} = C_{w_{j+1}}$ then
$g_j \notin D_{w_j}$ for any $j \in [1,k-1]$. Then, (\ref{eq:5}) can be characterized by a tuple
$$\{ |w_1(r_1)|, \ldots, |w_k(r_k)| \},$$
where $|w_j(r_j)| \in \mathbb{Z}^n,\ j \in [1,k]$. It is easy to see that any other representation
of $g$ as a reduced infinite word is very similar to (\ref{eq:5}). In particular, the following
statement holds.

\begin{lem}
\label{le:2.2.2}
Let $g \in G_{i+1} - G_i$ have two representations
$$g = g_1 \circ w_1(r_1) \circ g_2 \circ \cdots \circ w_k(r_k) \circ g_{k+1},$$
where $g_j \in G_i,\ j \in [1,k+1],\ w_j \in ST_i,\ j \in [1,k]$ and if $D_{w_j} = C_{w_{j+1}}$
then $g_j \notin D_{w_j}$ for any $j \in [1,k-1]$, and
$$g = h_1 \circ z_1(q_1) \circ h_2 \circ \cdots \circ z_m(q_m) \circ h_{m+1},$$
where $h_j \in G_i,\ j \in [1,m+1],\ z_j \in ST_i,\ j \in [1,m]$ and if $D_{z_j} = C_{z_{j+1}}$
then $h_j \notin D_{z_j}$ for any $j \in [1,m-1]$. Then $k = m$ and $w_j = z_j,\ j \in [1,k]$.
Moreover, $h_j = v_j \ast g_j \ast u_j$, where $v_j \in D_{w_{j-1}},\ u_j \in C_{w_j}$.
\end{lem}
{\it Proof.} We prove by induction on $k + m$. Observe that if $k + m = 0$ then there is nothing to
prove, so, let $k + m > 0$. We have
$$(g_1 \circ w_1(r_1) \circ \cdots \circ w_k(r_k) \circ g_{k+1}) \ast (h_1 \circ z_1(q_1) \circ
\cdots \circ z_m(q_m) \circ h_{m+1})^{-1} = \varepsilon,$$
that is,
$$(g_1 \circ w_1(r_1) \circ \cdots \circ w_k(r_k) \circ g_{k+1}) \ast (h_{m+1}^{-1} \circ
z_m^{-1}(-q_m) \circ \cdots \circ z_1^{-1}(-q_1) \circ h_1^{-1}) = \varepsilon,$$
which is only possible if $D_{w_k}^+ = D_{z_m}^+$ provided $d = g_{k+1} \ast h_{m+1}^{-1} \in
D_{w_k}$. By definition of $ST_i$ and by Theorem \ref{th:1} (6) we have $w_k = z_m$.
Let $c = w_k(r_k) \ast (g_{k+1} \ast h_{m+1}^{-1}) \ast z_m^{-1}(-q_m) \in C_{w_k}$, so we have
$$(g_1 \circ w_1(r_1) \circ \cdots \circ w_{k-1}(r_{k-1}) \circ g_k) \ast c$$
$$\ast (h_m^{-1} \circ z_{m-1}^{-1}(-q_{m-1}) \circ \cdots \circ z_1^{-1}(-q_1) \circ h_1^{-1}) =
\varepsilon.$$
Now, we have two representations of $h \in G_{i+1} - G_i$, where
$$h = g_1 \circ w_1(r_1) \circ \cdots \circ w_{k-1}(r'_{k-1}) \circ g'_k,$$
$$h = h_1 \circ z_1(q_1) \circ \cdots \circ z_{m-1}(q_{m-1}) \circ h_m,$$
where $w_{k-1}(r'_{k-1}) \circ g'_k = (w_{k-1}(r_{k-1}) \circ g_k) \ast d$, and by the induction
hypothesis we get $m - 1 = k - 1,\ w_j = z_j,\ j \in [1,k-1]$. So, the statement follows.

\hfill $\square$

From Lemma \ref{le:2.2.1} it follows that any representation of $g \in G_{i+1} - G_i$ can be associated
with a particular $k$-tuple for a fixed $k \in \mathbb{N}$. Now, representation (\ref{eq:5}) is called
the {\em normal form} if the $k$-tuple
$$\{ |w_1(r_1)|, \ldots, |w_k(r_k)| \}$$
is maximal in the left lexicographic order among $k$-tuples corresponding to all representations of
$g$. Existence of the normal representation (form) follows from Lemma 16 \cite{KMRS}. The (uniquely
defined) number $k$ is called the {\em level $i+1$ syllable length} of $g$. Elements of $G_i$ are said
to have level $i+1$ syllable length $0$. Normal forms satisfies the following properties.

\begin{rmk}
\label{re:2.2.1}
Normal forms are not unique since, while $w_j$ are defined uniquely by Lemma \ref{le:2.2.2},
elements $g_j$ are only defined up to multiplication $v \ast g_j \ast u$, where $v \in D(w_{j-1})$ and
$u \in C(w_j)$.
\end{rmk}

\begin{lem}
\label{le:2.2.3}
If
$$g = g_1 \circ w_1(r_1) \circ g_2 \circ \cdots \circ w_k(r_k) \circ g_{k+1}$$
is a normal form for $g$ then
\begin{enumerate}
\item[(a)] $g_1$ does not have $u_{w_1}^{\pm 1}$ as a terminal subword,
\item[(b)] $w_{i-1} \circ g_1$ does not have $u_{w_i}^{\pm 1}$ as a terminal subword for $i \in [2,k]$,
\item[(c)] $g_{i+1} \circ w_{i+1}(r_{i+1})$ does not have $u_{w_{i+1}}^{\pm 1}$ as an initial subword
for $i \in [1,k-1]$,
\item[(d)] $g_{k+1}$ does not have $v_{w_k}^{\pm 1}$ as an initial subword.
\end{enumerate}
\end{lem}
{\it Proof.} Since there is no cancellation in the products $g_i \circ w_i(r_i),\ w_i(r_i) \circ
g_{i+1},\ i \in [1,k]$, if at least one of the properties listed above does not hold then at least
one of $|w_i(r_i)|,\ i \in [1,k]$ can be made bigger by a finite sum of $|u|$ (even at the expense
of $|w_{i+1}(r_{i+1})|$), where $u$ can be taken to be a generator of any $C \in {\cal C}_i$ of
maximal height. This implies that $\{ |w_1(r_1)|, \ldots, |w_k(r_k)| \}$ is not maximal in the left
lexicographic order - a contradiction.

\hfill $\square$

Suppose
\begin{equation}
\label{eq:6}
g = g_1 \circ w_1(r_1) \circ g_2 \circ \cdots \circ w_k(r_k) \circ g_{k+1}
\end{equation}
is a normal form for $g$, where $g$ is cyclically reduced. Then, obviously
\begin{equation}
\label{eq:7}
\begin{array}{rl} g^2 =& g_1 \circ w_1(r_1) \circ g_2 \circ \cdots \circ w_k(r_k) \circ (g_{k+1}
\circ g_1)\\
& \phantom{a}\\
& \circ\ w_1(r_1) \circ g_2 \circ \cdots \circ w_k(r_k) \circ g_{k+1}
\end{array}
\end{equation}
is a representation of $g^2$ which may not be a normal form. So, we call (\ref{eq:6}) {\em cyclically
normal} if (\ref{eq:7}) is normal.

\begin{lem}
\label{le:2.2.4}
Let $u \in C^+$ and $v \in D^+$ be generators of $C, D \in {\cal C}_i,\ C \neq D$ of maximal height
such that $t^{-1} \ast u \ast t = v$ for some $t \in \langle T_i \rangle$. Let $w \in G_i$ be such
that $|w| < |u| = |v|$.
\begin{enumerate}
\item[(a)] If $u^\alpha \ast w \ast u^\beta = u^\alpha \circ w \circ u^\beta,\ \alpha, \beta \in
\mathbb{Z}$ and $w \circ u^\beta$ has $u^{\pm 1}$ as an initial subword then $u^\alpha \circ w$
cannot have $u^{\pm 1}$ as a terminal subword provided $w \notin C$.
\item[(b)] If $u^\alpha \ast w \ast v^\beta = u^\alpha \circ w \circ v^\beta,\ \alpha, \beta \in
\mathbb{Z}$ and $w \circ v^\beta$ has $u^{\pm 1}$ as an initial subword then $u^\alpha \circ w$
cannot have $v^{\pm 1}$ as a terminal subword.
\end{enumerate}
\end{lem}
{\it Proof.} {\bf (a)} Suppose $\beta > 0$.

If $w \circ u^\beta$ has $u$ as an initial segment then $\alpha > 0$ and $u = w \circ u_1,\ u = u_1
\circ u_2$. Observe that $u^\alpha \circ w$ cannot end with $u^{-1}$ but if it does end with $u$ then
$u = u_3 \circ w,\ u = u_4 \circ u_3$. Combining all four representations of $u$ we get $w = u_2$
which implies $u = w \circ u_1 = u_1 \circ w$, that is, $w \in C$ - a contradiction.

If $w \circ u^\beta$ has $u^{-1}$ as an initial segment then $\alpha < 0$ and $u^{-1} = w \circ u_1,\
u = u_1 \circ u_2$ which implies $u_1 = \varepsilon$ and $w = u_2^{-1} \in C$ - a contradiction.

\smallskip

Suppose $\beta < 0$.

If $w \circ u^\beta$ has $u$ as an initial segment then $\alpha > 0$ and $u = w \circ u_1,\ u^{-1} =
u_1 \circ u_2$ which implies $u_1 = \varepsilon$ and $w = u_2^{-1} \in C$ - a contradiction.

If $w \circ u^\beta$ has $u^{-1}$ as an initial segment then $\alpha < 0$ and $u^{-1} = w \circ u_1,\
u^{-1} = u_1 \circ u_2$. Observe that $u^\alpha \circ w$ cannot end with $u$ but if it does end with
$u^{-1}$ then $u^{-1} = u_3 \circ w,\ u^{-1} = u_4 \circ u_3$. Combining all four representations of
$u^{-1}$ we get $w = u_2$ which implies $u^{-1} = w \circ u_1 = u_1 \circ w$, that is, $w \in C$ -
a contradiction.

\smallskip

{\bf (b)} Suppose $\alpha > 0,\ \beta > 0$. If $w \circ v^\beta$ has $u$ as an initial segment then
$u = w \circ v_1,\ v = v_1 \circ v_2,\ |w| = |v_2|$. Observe that $u^\alpha \circ w$ cannot end with
$v^{-1}$ but if it does end with $v$ then $v = u_2 \circ w,\ u = u_1 \circ u_2,\ |w| = |u_1|$. Now
it follows that $w = u_1,\ v_1 = u_2$ and $v = u_2 \circ u_1$ which is a contradiction since $C$ and
$D$ are not conjugate in $G_i$.

\smallskip

Suppose $\alpha > 0,\ \beta < 0$. If $w \circ v^\beta$ has $u$ as an initial segment then $u = w
\circ v_1,\ v^{-1} = v_1 \circ v_2,\ |w| = |v_2|$. Observe that $u^\alpha \circ w$ cannot end with
$v$ but if it does end with $v^{-1}$ then $v^{-1} = u_2 \circ w,\ u = u_1 \circ u_2,\ |w| = |u_1|$.
Now it follows that $w = v_2,\ v_1 = u_2$ and $u = v_2 \circ v_1$, that is, $u$ is a conjugate of
$v^{-1}$. Finally, since $u$ is a conjugate of $v$ it follows that $v$ is a conjugate of $v^{-1}$ -
a contradiction.

\smallskip

The cases $\alpha < 0,\ \beta > 0$ and $\alpha < 0,\ \beta < 0$ can be considered in a similar way.

\hfill $\square$

\begin{lem}
\label{le:2.2.5}
Let
$$g = g_1 \circ w_1(r_1) \circ g_2 \circ \cdots \circ w_k(r_k) \circ g_{k+1}$$
be a normal form of a cyclically reduced $g \in G_{i+1} - G_i$. Then there exists a cyclic permutation
of $g$ whose normal form is cyclically normal unless $k = 1$ and $g_2 \ast g_1 \in D_{w_1} = C_{w_1}$.
\end{lem}
{\it Proof.} Without loss of generality we can assume $g_1 = \varepsilon$ and that $g_{k+1}$ does
not have $u_{w_1}$ as a terminal subword (using cyclic permutation we can always obtain these
properties).

If $v_{w_k} = u_{w_1}$ and $g_{k+1} \in D_{w_k}$ then applying a cyclic permutation we reduce the
number of entries of elements from $\langle T_i \rangle$ provided $k > 1$. Otherwise, $g = w_1(r_1)
\circ g_2,\ g_2 = g_{k+1} \in D_{w_1}$.

So we can assume that either $v_{w_k} \neq u_{w_1}$ or $g_{k+1} \notin D_{w_k}$

\begin{enumerate}
\item $|g_{k+1}| \geqslant |u_{w_1}|$

Since $w_1(r_1) \circ g_2 \circ \cdots \circ w_k(r_k) \circ g_{k+1}$ is normal it follows that
$g_{k+1}$ does not have $u_{w_k}$ as an initial and $u_{w_1}$ as a terminal subword. Hence,
$$w_1(r_1) \circ g_2 \circ \cdots \circ w_k(r_k) \circ g_{k+1} \circ w_1(r_1) \circ g_2 \circ \cdots
\circ w_k(r_k) \circ g_{k+1}$$
is normal.

\item $|g_{k+1}| < |u_{w_1}|$

If $g_{k+1} \circ w_1(r_1)$ does not have $v_{w_k}$ as an initial segment then
$$w_1(r_1) \circ g_2 \circ \cdots \circ w_k(r_k) \circ g_{k+1}$$
is cyclically normal.

Suppose $g_{k+1} \circ w_1(r_1)$ has $v_{w_k}$ as an initial segment. We take a cyclic
permutation
$$g' = g_2 \circ w_2(r_2) \circ \cdots \circ w_k(r_k) \circ g_{k+1} \circ w_1(r_1)$$
of $g$ so that
$$g_2 \circ w_2(r_2) \circ \cdots \circ w_k(r_k+1) \circ b \circ w_1(r_1-1)$$
is a normal form for $g'$, where $v_{w_k} = g_{k+1} \circ a,\ u_{w_1} = a \circ b$. But then
$$(g_2 \circ w_2(r_2) \circ \cdots \circ w_k(r_k+1) \circ b \circ w_1(r_1-1)) \circ (g_2 \circ
w_2(r_2) \circ$$
$$\cdots \circ w_k(r_k+1) \circ b \circ w_1(r_1-1))$$
is a normal form for $(g')^2$ by Lemma \ref{le:2.2.4}.
\end{enumerate}

\hfill $\square$

\subsection{Normal forms in $G$}
\label{sub:2.3}

Recall that $G$ can be represented as a union of a finite series of groups~(\ref{eq:2})
$$G_1 < G_2 < \cdots < G_{K(G)} = G$$
satisfying the properties shown in Theorem \ref{th:1}. In particular $G_1 = \langle X \rangle$ is a
free group.

Using the series (\ref{eq:2}) we can introduce a {\em standard decomposition} $\pi(g)$ for $g \in
G$, where $\pi(g) \in ({\cal B}(G)^{\pm 1})^*$ and
$${\cal B}(G) = X \bigcup \left(\bigcup_i T_i \right) \bigcup \left( \{ u^n \mid n \in \mathbb{Z},\
u \in U \}\right),$$
where
$$U = \bigcup_i U_i,$$
$$U_i = \{ u \in C \mid C \in {\cal C}_i,\ u\ {\rm is\ a\ generator\ of}\ C\ {\rm of\ maximal\
height} \}.$$
If $g \in G_1$ then
$$g = x_1 \cdots x_k,\ x_i \in X^{\pm 1},$$
which a reduced word in a free group $G_1$, so we put
$$\pi(g) = x_1 \cdots x_k.$$
Assuming inductively that $\pi(f)$ is constructed for any $f \in G_i$ we consider $g \in G_{i+1} -
G_i$. By the results of the previous subsection we can compute a normal form for
$g$
$$g = g_1 \circ w_1(r_1) \circ g_2 \circ \cdots \circ w_k(r_k) \circ g_{k+1}$$
satisfying the properties shown in Lemma \ref{le:2.2.3}. Since $g_i,\ i \in [1,k+1]$ are elements of
$G_i$ it follows that $\pi(g_i),\ i \in [1,k+1]$ are defined and we put
$$\pi(g) = \pi(g_1)\ w_1(r_1)\ \pi(g_2)\ \cdots\ w_k(r_k)\ \pi(g_{k+1}).$$

\smallskip

Observe that in general, for a cyclically reduced $g \in G$ we have $\pi(g^2) \neq \pi(g)^2$. In
particular, if $u \in G_{i+1}$ is such that there exists $t \in T_{i+1}$ such that $u = u_t$ then
$\pi(u^2)$ does not have to be equal to $\pi(u)^2$. On the other hand, by Lemma \ref{le:2.2.5}
every $g \in G_{i+1} - G_i$ has a cyclic permutation whose normal form is cyclically normal unless
$g = w \circ h$, where $w \in ST_i$ and $h \in D_w = C_w$. Since all centralizers from
${\cal C}_{i+1}$ can be taken up to conjugation (cyclic permutation) then we can always assume that
$\pi(u^2) = \pi(u)^2$ unless $u = w \circ h$, where $w \in ST_i$ and $h \in D_w = C_w$. Observe that
in the latter case $t$ has any natural power of $u$ as an initial and terminal subword and without
loss of generality we can assume $u = w$ which is the generator of maximal height in $C_t =
C_{G_{i+1}}$. Hence, $\pi(w^2) = \pi(w)^2$.

\section{${\cal B}(G)$-graphs}
\label{sec:3}

Let us fix a group $G$ with a free regular length function in $\mathbb{Z}^n$ for the rest of this
section. According to the previous section, every $g \in G$ can be written as a word in the alphabet
${\cal B}(G)^{\pm 1}$. In this section we introduce graphs whose edges are labeled by ${\cal B}(G)$
and show how such graphs can be associated with subgroups of $G$.

\subsection{Labeled graphs}
\label{subs:3.1}

Using the notation introduced in the previous section we adjust basic notions from \cite{KM} to the
case of groups acting freely on $\mathbb{Z}^n$-trees.

\begin{defn}
\label{de:3.1.1}
By a ${\cal B}(G)$-labeled directed graph ({\em ${\cal B}(G)$-graph}) $\Gamma$ we mean the following.
\begin{enumerate}
\item [(1)] $\Gamma$ is a combinatorial graph, where every edge has a direction and is labeled by a
letter from ${\cal B}(G)$, denoted $\mu(e)$,
\item[(2)] for each edge $e$ of $\Gamma$ we denote the origin of $e$ by $o(e)$ and the terminus of
$e$ by $t(e)$.
\end{enumerate}
\end{defn}

For each edge $e$ of ${\cal B}(G)$-graph we introduce a formal inverse $e^{-1}$ of $e$ with label
$\mu(e)^{-1}$ and the endpoints defined as $o(e^{-1}) = t(e),\ t(e^{-1}) = o(e)$, that is, the
direction of $e^{-1}$ is reversed with respect to the direction of $e$. For the new edge $e^{-1}$ we
set $(e^{-1})^{-1} = e$. The new graph, endowed with this additional structure we denote by
$\widehat{\Gamma}$. In fact, usually we will abuse notation by disregarding the difference between
$\Gamma$ and $\widehat{\Gamma}$.

Now we have a partition $E(\widehat{\Gamma}) = E(\Gamma) \cup \overline{E(\Gamma)}$ and we say that
edges of $\Gamma$ are {\em positively oriented} in $\widehat{\Gamma}$, while their formal inverses
$e^{-1}$ are {\em negatively oriented} in $\widehat{\Gamma}$.

\begin{defn}
\label{de:3.1.2}
A path $p$ in $\Gamma$ is a sequence of edges $p = e_1 \cdots e_k$, where each $e_i$ is an edge of
$\widehat{\Gamma}$ and the origin of each $e_i$ is the terminus of $e_{i-1}$.
\end{defn}

Observe that $\mu(p) = \mu(e_1) \cdots \mu(e_k)$ is a word in the alphabet ${\cal B}(G)$ and we
denote by $\overline{\mu(p)}$ the reduced infinite word $\mu(e_1) \ast \cdots \ast \mu(e_k)$.

\smallskip

We will be using two different notions of the length of a path $p = e_1 \cdots e_k$ in $\Gamma$
\begin{enumerate}
\item {\em combinatorial length} $|p|$ set equal to $k$, and
\item {\em word length} $wl(p) = \sum_{i = 1}^k l(\mu(e_i))$.
\end{enumerate}

In fact, from these two definitions above two possible meanings of irreducible path arise: irreducible
path in combinatorial sense and irreducible path in the sense that its label is viewed as an infinite
word in $CDR(\mathbb{Z}^n,X)$ is reduced. Here are formal definitions.

\begin{defn}
\label{de:3.1.3}
A path $p = e_1 \cdots e_k$ in a ${\cal B}(G)$-graph $\Gamma$ is called reduced if $e_i \neq
e_{i+1}^{-1}$ for all $i \in [1,k-1]$.
\end{defn}

\begin{defn}
\label{de:3.1.4}
A path $p = e_1 \cdots e_k$ in a ${\cal B}(G)$-graph $\Gamma$ is called {\em label reduced} if
\begin{enumerate}
\item $p$ is reduced,
\item if $t = \mu(e_l) = \mu(e_m)^{-1},\ t \in T_i^{\pm 1},\ l < m,\ l,m \in [1,k]$ and $q = e_{l+1}
\cdots e_{m-1}$ does not contain edges labeled by $T_j^{\pm 1},\ j \geqslant i$ then $w = \overline{\mu(q)}
\notin D_t$, moreover, $t \ast w = t \circ w,\ w \ast t^{-1} = w \circ t^{-1}$.
\end{enumerate}
\end{defn}

Recall that a graph labeled by letters from $X^{\pm}$ defines a language of words over $X^\pm$. This
language can be put into correspondence with a subgroup of a free group $F(X)$. In the present
subsection we generalize this concept to ${\cal B}(G)$-graphs.

\begin{defn}
\label{de:3.1.5}
Let $\Gamma$ be a ${\cal B}(G)$-graph and let $v,v'$ be vertices of $\Gamma$. We define the language
of $\Gamma$ with respect to $v$ to be
$$L(\Gamma, v) = \{\overline{\mu(p)} \mid\ p\ {\rm is\ a\ reduced\ path\ in}\ \Gamma\ {\rm from}\ v\
{\rm to}\ v\}.$$
We also put
$$L(\Gamma, v, v') = \{\overline{\mu(p)} \mid\ p\ {\rm is\ a\ reduced\ path\ in}\ \Gamma\ {\rm from}\
v\ {\rm to}\ v'\}.$$
\end{defn}

If $w$ belongs to $L(\Gamma, v)$, we will also sometimes say that $w$ is {\em accepted} by
$(\Gamma, v)$ (or just by $\Gamma$ if $v$ is fixed).

The following result establishes a connection between ${\cal B}(G)$-graphs and subgroups of $G$.

\begin{lem}
\label{le:3.1.1}
Let $\Gamma$ be a finite ${\cal B}(G)$-graph and let $v \in V(\Gamma)$. Then $L(\Gamma, v)$ is a
subgroup of $G$.
\end{lem}
{\it Proof.} Straightforward verification.

\hfill $\square$

Similarly, $L(\Gamma, v, v_1)$ is a left $L(\Gamma,v')$-coset and a right $L(\Gamma,v)$-coset of
$G$.

\subsection{Free foldings}
\label{subs:3.2}

Here we define free (partial) foldings and partially folded $(\mathbb{Z}^n,X)$-graphs. Observe that
the definition of a partial folding below is exactly the same as the corresponding definition of a
folding in free groups (see \cite{KM}).

Let $\Gamma$ be a ${\cal B}(G)$-graph. Suppose $v_0$ is a vertex of $\Gamma$ and $f_1,\ f_2$ are two
distinct edges of $\widehat{\Gamma}$ such that $o(f_1) = o(f_2) = v_0,\ \mu(f_1) = \mu(f_2) \in
{\cal B}^\pm$. Let $h_i$ be the positive edge of $\Gamma$ corresponding to $f_i$ (that is, $h_i =
f_i$ if $f_i$ is positive and $h_i = {f_i}^{-1}$ if $f_i$ is negative).

Let $\Delta$ be a ${\cal B}(G)$-graph with the following sets of vertices and edges.
$$V(\Delta) = (V(\Gamma) - \{t(f_1),t(f_2)\}) \cup \{v\},\ \ E(\Delta) = (E(\Gamma) - \{h_1,h_2\})
\cup \{h\}.$$
The endpoints and arrows for the edges of $\Delta$ are defined in the following way. Let $e \in
E(\Delta),\ e \neq h$ then
\begin{enumerate}
\item we put $o_{\Delta}(e) = o_{\Gamma}(e)$ if $o_{\Gamma}(e) \ne t(f_i)$ and  $o_{\Delta}(e) = v$
if $o_{\Gamma}(e) =  t(f_i)$ for some $i$,
\item we put $t_{\Delta}(e) = t_{\Gamma}(e)$ if $t_{\Gamma}(e)\ne t(f_i)$ and $t_{\Delta}(e) = v$ if
$t_{\Gamma}(e) =  t(f_i)$ for some $i$.
\end{enumerate}
For the edge $h$ we put $o_{\Delta}(h) = v_0,\ t_{\Delta}(h) = v$ if $h_1 = f_1,\ h_2 = f_2$ and
$o_{\Delta}(h) = v,\ t_{\Delta}(h) = v_0$ otherwise.

We define labels on the edges of $\Delta$ as follows: $\mu_{\Delta}(e) = \mu_{\Gamma}(e)$ if $e \ne
h$ and $\mu_{\Delta}(h) = \mu_{\Gamma}(h_1) = \mu_{\Gamma}(h_2)$.

In other words we obtain $\Delta$ by identification of two edges $f_1$ and $f_2$ in $\Gamma$. In
this situation we say that $\Delta$ is obtained from $\Gamma$ by a {\em free folding} (or by {\em
freely folding the edges $f_1$ and $f_2$}).

There can be introduced a notion of a {\em morphism} between two ${\cal B}(G)$-graphs. That is, if
$\Gamma_1,\ \Gamma_2$ are ${\cal B}(G)$-graphs then a map $\theta : \Gamma_1 \rightarrow \Gamma_2$
is called a {\em morphism} of ${\cal B}(G)$-graphs, if $\theta$ sends vertices to vertices, directed
edges to directed edges, preserves labels of directed edges, and has the property that $o(\theta(e))
= \theta(o(e)),\ t(\theta(e)) = \theta(t(e))$ for any edge $e$ of $\Gamma_1$.

If $\phi$ is a free folding defined above then it is easy to see that  $\phi$ is a morphism between
$\Gamma$ and $\Delta$.

\begin{lem}
\label{le:3.2.1}
Let $\Gamma_1$ be a ${\cal B}(G)$-graph obtained by a free folding from a graph $\Gamma$. Let $v, v'$
vertices of $\Gamma$ and $v_1,v'_1$ be the corresponding vertices of $\Gamma_1$. Then the following
hold.
\begin{enumerate}
\item If $\Gamma$ is connected then $\Gamma_1$ is connected.
\item Let $p$ be the path from $v$ to $v'$ in $\Gamma$ with label $w$. Then the edgewise image of
$p$ in $\Gamma_1$ is a path from $v_1$ to $v'_1$ with label $w$.
\item If $\Gamma$ is a finite ${\cal B}(G)$-graph, then the number of edges in $\Gamma_1$ is one
less than the number of edges in $\Gamma$, that is, any free folding decreases the number of edges
in $\Gamma$.
\end{enumerate}
\end{lem}
{\it Proof.} Follows directly from the definition of a free folding.

\hfill $\square$

\begin{defn}
\label{de:3.2.2}
${\cal B}(G)$-graph $\Gamma$ is called {\em freely folded} if there exist no two edges $e_1$ and
$e_2$ in $\Gamma$ with $\mu(e_1) = \mu(e_2)$ such that $o(e_1) = o(e_2)$ or $t(e_1) = t(e_2)$.
\end{defn}

Obviously, $\Gamma$ is a freely folded ${\cal B}(G)$-graph if and only if one cannot perform any
free folding in $\Gamma$. Moreover the following proposition is true.

\begin{lem}
\label{pr:3.2.1}
Let $\Gamma$ be a ${\cal B}(G)$-graph, which has only a finite number of edges. Then there exists a
freely folded ${\cal B}(G)$-graph $\Delta$, which can be obtained from $\Gamma$ by a finite number
of free foldings.
\end{lem}
{\it Proof.} Since $\Gamma$ has a finite number of edges by Lemma \ref{le:3.2.1}, any
${\cal B}(G)$-graph $\Gamma_1$ obtained from $\Gamma$ by a free folding has fewer edges. This
provides one with an inductive argument based on the number of edges in $\Gamma$.

\hfill $\square$

\begin{lem}
\label{le:3.2.2}
Let $\Gamma$ be a finite ${\cal B}(G)$-graph and let $v,v' \in V(\Gamma)$. Let $\Delta$ be a
${\cal B}(G)$-graph obtained from $\Gamma$ by a single free folding so that $v_1,v'_1 \in V(\Delta)$
correspond to $v,v'$. Then
$$L(\Gamma,v,v') = L(\Delta,v_1,v_1').$$
\end{lem}
{\it Proof.} Similar to the proof of Lemma 3.4 \cite{KM}.

\hfill $\square$

\subsection{$C$-components}
\label{subs:3.3}

In the present subsection we fix $C \in {\cal C}_m$ and concentrate on subgraphs of
${\cal B}(G)$-graphs which consist of edges labeled by exponents $u^k,\ u \in U_i \cap C,\ i \leqslant m$.

\smallskip

Let $\Gamma$ be a ${\cal B}(G)$-graph and let $C \in {\cal C}_m$. We say that $v_0,v_1 \in V(\Gamma)$
are {\em $C$-equivalent} (denote $v_0 \sim_C v_1$) if there exists a path $p$ in $\Gamma$ such that
$o(p) = v_0,\ t(p) = v_1$, and every edge of $p$ is labeled by a power of $u \in U_i \cap C,\ i \leqslant
m$. In other words, if $v_0 \sim_C v_1$ then there exists a path connecting them, whose reduced
label is an element of $C$.

One can take the subgraph of $\Gamma$ spanned by vertices $C$-equivalent to $v \in V(\Gamma)$ and
remove from it all edges labeled by anything except for $u \in U_i \cap C,\ i \leqslant m$. The resulting
subgraph of $\Gamma$ we denote by $Comp_C(v)$ and call the {\em $C$-component of $\ v$}.

\begin{defn}
\label{de:3.3.1}
Let $\Gamma$ be a ${\cal B}(G)$-graph and $v \in V(\Gamma),\ v_0,v_1 \in V(Comp_C(v))$. We define a
set $H_C(v_0)$ associated with $v_0$ as
$$H_C(v_0) = \{\overline{\mu(p)} \mid p\ {\rm is\ a\ reduced\ path\ in}\ Comp_C(v)\ {\rm from}\
v_0\ {\rm to}\ v_0 \}.$$
We also put
$$H_C(v_0,v_1) = \{\overline{\mu(p)} \mid p\ {\rm is\ a\ reduced\ path\ in}\ Comp_C(v)\ {\rm from}\
v_0\ {\rm to}\ v_1 \}.$$
\end{defn}

Observe that even when $p$ is a reduced path in $Comp_C(v)$ its label $\overline{\mu(p)}$ may be the
empty infinite word.

\begin{lem}
\label{le:3.3.1}
Let $\Gamma$ be a ${\cal B}(G)$-graph and $v \in V(\Gamma),\ v_0 \in V(Comp_C(v))$. Then
\begin{enumerate}
\item [(1)] $H_C(v_0)$ is a subgroup of $C$,
\item [(2)] if $v_1 \in V(Comp_C(v))$ then $H_C(v_0) = H_C(v_1)$.
\end{enumerate}
\end{lem}
{\it Proof.} Straightforward verification.

\hfill $\square$

It follows from Lemma \ref{le:3.3.1} that one can associate a subgroup of $C$ with any finite
$C$-component $Q$ in a ${\cal B}(G)$-graph $\Gamma$. We can denote this subgroup by $H_C(Q)$.

\subsection{Reduced $C$-components}
\label{subs:3.4}

Suppose $Q$ is a partially folded $C$-component of a ${\cal B}(G)$-graph $\Gamma$ and fix $v_0 \in
V(Q)$. For a vertex $v \in V(Q)$ take a path $p_v$ such that $o(p_v) = v_0,\ t(p_v) = v$. So
$\overline{\mu(p_v)} \in C$ and obviously for any other path $q_v$ such that $o(q_v) = v_0,\ t(q_v)
= v$ we have $\overline{\mu(p_v)} \ast \overline{\mu(q_v)} \in H_C(Q)$ since $p_v q_v^{-1}$ is
a loop at $v_0$. Hence, every vertex $v \in V(Q)$ can be associated with a coset $H_C(v_0,v) =
\overline{\mu(p_v)} \ast H_C(Q)$ in $C$ by $H_C(Q)$. Denote this system of cosets by $C_Q(v_0)$.

It can happen that two distinct vertices of $Q$ correspond to the same coset in $C_Q(v_0)$, that is,
$$\overline{\mu(p_{v_1})} \ast H_C(Q) = \overline{\mu(p_{v_2})} \ast H_C(Q),\ \ v_1 \neq v_2.$$
Observe that in this case $v_1$ and $v_2$ correspond to the same coset in $C_Q(v)$ for any $v \in
V(Q)$ and we can apply a {\em simple reduction} of $Q$, which is an identification of the vertices
$v_1$ and $v_2$. It is easy to see that simple reductions are morphisms of graphs.

\begin{lem}
\label{le:3.4.1}
Let $Q$ be a $C$-component of a finite ${\cal B}(G)$-graph $\Gamma$. Let $v_1, v_2 \in V(Q)$
correspond to the same coset in $C_Q(v_0)$ for some $v_0 \in V(Q)$, and let $\Delta$ be a
${\cal B}(G)$-graph obtained from $\Gamma$ by a simple reduction $\phi$ of $Q$ such that $\phi(v_1)
= \phi(v_2)$. If $v_0, v \in V(\Gamma)$ and $v'_0=\phi(v_0),v' = \phi(v)$ then
$$L(\Gamma,v_0,v) = L(\Delta,v'_0,v').$$
\end{lem}
{\it Proof.} Straightforward verification similar to Lemma \ref{le:3.3.1}.

\hfill $\square$

If every coset in $C_Q(v_0),\ v_0 \in V(Q)$ corresponds to a unique vertex of $Q$ then we say that
$Q$ is a {\em reduced} $C$-component. Since simple reductions decrease the number of vertices it
follows that after finitely many reductions we obtain a reduced $C$-component.

\subsection{$ST_m$-complexes}
\label{subs:3.5}

A path $p$ in a ${\cal B}(G)$-graph $\Gamma$ we call an {\em $ST_m$-path} if $\mu(p) \in ST_m$, that
is, if $p = e_1 \cdots e_k$ then $\mu(e_i) = t_i \in T_m^{\pm 1}$ and $D_{t_i} = C_{t_{i+1}},\ i
\in [1,k-1]$. If $Q$ is a non-trivial $C$-component of $\Gamma$, where $C \in {\cal C}_m$, then we
call $Q$ {\em attached to an $ST_m$-path $p$} if $o(p)$ belongs to $V(Q)$ and $C = C_w$, where $w =
\mu(p) \in ST_m$. A connected subgraph $P$ of $\Gamma$ is called an {\em $ST_m$-subgraph} if for
every two vertices $v_1, v_2 \in V(P)$ there exists a path
$$p = s_0 p_1 s_1 p_2 \cdots p_k s_k,$$
where $o(p) = v_1,\ t(p) = v_2$, every $p_i$ is a non-trivial $ST_m$-path, every $s_i$ is a path in a
$C_i$-component $Q_i$ of $\Gamma$, and $Q_0$ is attached to $p_1$, $Q_i$ is attached both to $p_i^{-1}$
and $p_{i+1}$ for $i \in [1,k-1]$, and $Q_k$ is attached to $p_k^{-1}$. It is easy to see that every
$ST_m$-path is an $ST_m$-subgraph which does not contain
non-trivial ${\cal C}_m$-components. Moreover, if $\Gamma$ is freely folded then every $ST_m$-subgraph
without non-trivial ${\cal C}_m$-components is an $ST_m$-path.

A maximal (by inclusion) $ST_m$-subgraph we call an {\em $ST_m$-complex}.

\smallskip

Let $Q_1$ be a $C_1$-component and $Q_2$ a $C_2$-component of an $ST_m$-complex $P$. Let $p$ be any
path from $Q_1$ to $Q_2$ in $P$. For any loop $q$ at $t(p)$ in $Q_2$ we have $c_1 = \overline{\mu(p)}
\ast \overline{\mu(q)} \ast \overline{\mu(p)}^{-1} \in C_1$ but it may happen that $c_1 \notin
H_{C_1}(Q_1)$. In other words, at a vertex of $Q_1$ there may be a loop in $P$ whose reduced label
belongs to $C_1$ but is not readable in $Q_1$. In this case we can enlarge $H_{C_1}(Q_1)$ by adding
a loop at $o(p)$ labeled by $c_1$. More precisely, let $q = e_1 \cdots e_k$, where $\mu(e_i) =
u_i^{n_i},\ u_i \in U_m \cap C_2$ for each $i \in [1,k]$. We form
a loop $q' = f_1 \cdots f_k$ such that $\mu(f_i) = v_i^{n_i},\ v_i \in U_m \cap C_2$, where $v_i =
\overline{\mu(q)} \ast u_i \ast \overline{\mu(q)}^{-1}$ for $i \in [1,k]$, and attach it to $o(p)$.
Such an operation we call an {\em elementary loop translation}. That is, we translated a loop $q$
from $Q_2$ to $Q_1$. It is easy to see that an elementary loop translation is a graph morphism since
the resulting graph is the original one together with a loop attached to one of its vertices.

A loop translation can be performed even if $Q_1$ is a trivial ${\cal C}_m$-component of $P$, that
is, $Q_1$ is a single vertex.

\begin{lem}
\label{le:3.5.1}
Let $\Gamma$ and $\Delta$ be ${\cal B}(G)$-graphs such that $\Delta$ is obtained from $\Gamma$ by an
elementary loop translation $\phi$ in an $ST_m$-complex $P$ of $\Gamma$. If $v_0, v \in V(\Gamma)$
and $v'_0 = \phi(v_0), v' = \phi(v)$ then
$$L(\Gamma, v_0, v) = L(\Delta, v'_0, v').$$
\end{lem}
{\it Proof.} We prove this lemma for the case $v = v_0$. The argument can be easily adjusted to
accommodate the case $v \neq v_0$.

Let $p$ be a path from $Q_1$ to $Q_2$ in $\Gamma$, where $Q_1$ is a $C_1$-component and $Q_2$ a
$C_2$-component of the $ST_m$-complex $P$, and $\Delta$ is obtained from $\Gamma$ by translation of
a loop $q$ from $t(p) \in V(Q_2)$ to $o(p) \in V(Q_1)$. The new loop in $Q_1$ we denote by $q'$.
With abuse of notation we denote $\phi(p)$ and $\phi(q)$ again respectively by $p$ and $q$ in
$\Delta$.

Obviously $L(\Gamma,v) \subseteq L(\Delta,v')$ since $\Gamma$ is graph-isomorphic to $\phi(\Gamma)$.

Now let $r$ be a reduced loop at $v'$ in $\Delta$. If $r$ does not contain $q'$ then $\phi^{-1}(r)$
is a loop at $v$ in $\Gamma$ and $\mu(r) = \mu(\phi^{-1}(r))$, that is, $\overline{\mu(r)} \in
L(\Delta,v')$. Now, suppose $r$ contains loops $q'$, so we subdivide it as follows
$$r = r_1 q' r_2 \cdots r_{k-1} q' r_k,$$
where $t(r_i) = o(r_{i+1}) = o(p),\ i \in [1,k-1]$. Hence, $r' = r_1 r_2 \cdots r_k$ is a loop at
$v'$ which does not contain $q'$ and hence there exists its preimage $\phi^{-1}(r') = \phi^{-1}(r_1)
\cdots \phi^{-1}(r_k)$ in $\Gamma$ which is a loop at $v$. Now, for each entry of $q'$ in $r$ we
insert a loop $p\ q\ p^{-1}$ in $P$ into $\phi^{-1}(r')$, that is, we consider
$$r'' = \phi^{-1}(r_1)\ (p\ q\ p^{-1})\ \phi^{-1}(r_2) \cdots \phi^{-1}(r_{k-1})\ (p\ q\ p^{-1})\
\phi^{-1}(r_k).$$
Since $\overline{\mu(q')} = \overline{p\ q\ p^{-1}}$ and $\overline{\mu(r_1)} =
\overline{\mu(\phi^{-1}(r_1))}$ then $\overline{\mu(r)} = \overline{\mu(r'')} \in L(\Gamma,v)$. That
is, $L(\Delta,v') \subseteq L(\Gamma,v)$.

\hfill $\square$

Since $H_{C_2}(Q_2)$ is finitely generated then it is enough to perform only finitely many elementary
loop translations to obtain
$$\overline{\mu(p)} \ast H_{C_2}(Q_2) \ast \overline{\mu(p)}^{-1} \subseteq H_{C_1}(Q_1)$$
in the resulting graph.

If an $ST_m$-complex $P$ contains finitely many ${\cal C}_m$-components (including trivial ones) then
after application of finitely many elementary loop translations we get
$$\overline{\mu(p_{i,j})} \ast H_{C_j}(Q_j) \ast \overline{\mu(p_{i,j})}^{-1} = H_{C_i}(Q_i),$$
where $C_i, C_j$ are (possibly trivial) ${\cal C}_m$-components of $P$ and $p_{i,j}$ is a path from
$Q_i$ to $Q_j$ in $P$. In this event we call $P$ a {\em balanced} $ST_m$-complex.

\smallskip

We call an $ST_m$-complex $P$ {\em reduced} if for every $C$-component $Q$ of $P$, $C \in {\cal C}_m$,
there are no two distinct edges $e_1, e_2 \in E(P)$ such that $o(e_1), o(e_2) \in V(Q),\ \mu(e_1) =
\mu(e_2) = t \in T_m^{\pm 1}$. Observe that the condition $e_1, e_2 \in E(P)$ implies $C = C_t$. It
is easy to see that the property ``to be reduced'' for an $ST_m$-complex is similar to the property
``to be freely folded'' for a ${\cal B}(G)$-graph defined above. The latter property can be obtained
by a sequence of free foldings, and now we introduce a similar operation, application of which
transforms an arbitrary $ST_m$-complex into a reduced one.

\smallskip

Let $Q$ be a $C$-component of an $ST_m$-complex $P$ and suppose there exist distinct edges $e, f$
such that $o(e), o(f) \in V(Q),\ \mu(e) = \mu(f) = t \in T_m^{\pm 1}$. Let $p = e_1 \cdots e_k$ be
a path in $Q$ such that $o(p) = o(e),\ t(p) = o(f),\ \mu(e_i) = u_i^{n_i}$, where $u_i \in U_m \cap
C$ is a power of a generator of $C$ for each $i \in [1,k]$. Observe that $C_t = C$. Now, we form a
path $q = f_1 \cdots f_k$ such that $\mu(f_i) = v_i^{n_i},\ v_i \in U_m \cap D_t$, where $v_i =
t^{-1} \ast u_i \ast t$ for $i \in [1,k]$. Finally, we remove the edge $e$ and attach $q$ to $\Gamma$
so that $o(q) = t(e),\ t(q) = t(f)$. The resulting graph we denote by $\Delta$ and say that $\Delta$
is obtained from $\Gamma$ by an {\em $ST_m$-folding}.

\begin{lem}
\label{le:3.5.2}
Let $\Gamma$ and $\Delta$ be ${\cal B}(G)$-graphs such that $\Delta$ is obtained from $\Gamma$ by an
$ST_m$-folding in an $ST_m$-complex $P$ of $\Gamma$. If $v_0,v \in V(\Gamma)$ and $v_0',v' \in
V(\Delta)$ correspond to $v_0,v$ respectively, then
$$L(\Gamma,v_0,v) = L(\Delta,v_0',v').$$
\end{lem}
{\it Proof.} We prove this lemma for the case $v = v_0$. The argument can be easily changed to
accommodate the case $v \neq v_0$.

We are going to use the notation from the definition of an $ST_m$-folding. That is, suppose there
exist distinct edges $e, f$ such that $o(e), o(f) \in V(Q),$ $\ \mu(e) = \mu(f) = t \in T_m^{\pm 1}$,
where $Q$ is a $C$-component of an $ST_m$-complex $P$ of $\Gamma$. Suppose $p$ is a path in $Q$ from
$o(e)$ to $o(f)$, whose edges are labeled by powers of generators of $C_t = C$, $q$ is a corresponding
path whose edges are labeled by powers of generators of $D_t$, and $\Delta$ is obtained from $\Gamma$
by deleting $e$ and attaching $q$.

With abuse of notation we are going to denote edges and vertices of $\Gamma$ and $\Delta$ which
correspond to each other by the same letters.

\smallskip

Let $r$ be a reduced loop at $v$ in $\Gamma$. We subdivide $r$ as follows
$$r = p_1 e^{\delta_1} p_2 \cdots p_{k-1} e^{\delta_{k-1}} p_k,$$
where $o(p_1) = t(p_k) = v$ and $\delta_i \in \{-1,1\},\ i \in [1,k-1]$. Since in $\Delta$ the path
$p f q^{-1}$ begins at $o(e)$ and ends at $t(e)$ then we can consider a loop $r'$ at $v'$ in
$\Delta$, where
$$r' = p_1 (p f q^{-1})^{\delta_1} p_2 \cdots p_{k-1} (p f q^{-1})^{\delta_{k-1}} p_k.$$
Now, since $\overline{\mu(p f q^{-1})} = \mu(e)$, we have $\overline{\mu(r)} = \overline{\mu(r')}$
and $L(\Gamma,v) \subseteq L(\Delta,v')$.

\smallskip

Now, let $r$ be a reduced loop at $v'$ in $\Delta$. Observe that since $q$ is connected to $\Delta$
only at the endpoints and $r$ is reduced it follows that if $r$ contains an edge (or its inverse)
which belongs to $q$ then $r$ contains the whole path $q$. Hence, we can subdivide $r$ as follows
$$r = p_1 q^{\delta_1} p_2 \cdots p_{k-1} e^{\delta_{k-1}} p_k,$$
where $o(p_1) = t(p_k) = v$ and $\delta_i \in \{-1,1\},\ i \in [1,k-1]$. Observe that the path
$e^{-1} p f$ begins at $o(q)$ and ends at $t(q)$, so we can consider a loop $r'$ at $v$ in $\Gamma$,
where
$$r' = p_1 (e^{-1} p f)^{\delta_1} p_2 \cdots p_{k-1} (e^{-1} p f)^{\delta_{k-1}} p_k.$$
Now, since $\overline{\mu(e^{-1} p f)} = \mu(q)$, we have $\overline{\mu(r)} = \overline{\mu(r')}$
and $L(\Delta,v') \subseteq L(\Gamma,v)$.

\hfill $\square$

\begin{rmk}{\rm Note that while $ST_m$-folding is not, strictly speaking, always a graph morphism,
it is a {\em generalized morphism of ${\cal B}(G)$-graphs,} that is a map $\theta : \Gamma_1
\rightarrow \Gamma_2$ between ${\cal B}(G)$-graphs $\Gamma_1,\ \Gamma_2$ such that $\theta$ sends
vertices to vertices, directed edges to directed paths, preserves labels of directed edges in the
sense that $\overline{\mu(\theta(e))}=\overline{\mu(e)}$, and has the property that $o(\theta(e)) =
\theta(o(e)),\ t(\theta(e)) = \theta(t(e))$ for any edge $e$ of $\Gamma_1$.}
\end{rmk}

If $P$ is a finite $ST_m$-complex of a ${\cal B}(G)$-graph $\Gamma$ then there exists a finite
sequence of $ST_m$-foldings transforming $P$ into a reduced $ST_m$-complex. Indeed, an $ST_m$-folding
reduces the number of edges of $P$ labeled by $T_m^{\pm 1}$, so the process stops in finitely many
steps.

Next, due to the properties of $T_m$, for every $C \in {\cal C}_m$ there are at most two elements
$t, s \in T_m$ attached to $C$. Hence, from the definition of $ST_m$-complex it follows that if $Q$
is a $C$-component of $P$ then it can be attached to edges labeled at most by two letters from
$T_m^{\pm 1}$. So, in the corresponding reduced $ST_m$-complex there can be at most two edges
outgoing from $Q$. Thus, if we contract all ${\cal C}_m$-components of $P$ into vertices then the
resulting graph $P / \sim_{{\cal C}_m}$ is a simple path (loop). This is summarized in the following
lemma.

\begin{lem}
\label{le:3.5.3}
Let $P$ be a balanced reduced $ST_m$-complex of a finite ${\cal B}(G)$-graph $\Gamma$. Then $P$ is a
simple path (loop) modulo ${\cal C}_m$-equivalence. That is, every ${\cal C}_m$-component of $P$ is
attached to at most two edges labeled by $T_m^{\pm 1}$. Moreover, every vertex of $P$ belongs to
some ${\cal C}_m$-component of $P$ and the subgroups associated to all ${\cal C}_m$-components are
isomorphic.
\end{lem}
{\it Proof.} Follows from the discussion above.

\hfill $\square$

\begin{lem}
\label{le:3.5.4}
Let $P$ be a balanced reduced $ST_m$-complex of a finite ${\cal B}(G)$-graph $\Gamma$. If $v \in V(P),
\ w \in ST_m,\ u\in C_w$ then there exists an algorithm which effectively determines if there
exists a path $p$ in $P$ such that $o(p) = v$ and $\overline{\mu(p)} = u \ast w$.
\end{lem}
{\it Proof.} If such a path $p$ exists, the element $w$ uniquely determines the sequence of
$ST_m$-edges $e_1, e_2, \ldots, e_l$ that $p$ traverses and, since $P$ is reduced, which $C$-component
$Q$ of $P$ contains the end-point $t(p)$. Let $v_1, v_2, \ldots, v_k$ be vertices of $Q$. For each
$v_i$, fix an arbitrary $ST_m$-path $q_i$ that originates at $v$, traverses $ST_m$-edges $e_1, e_2,
\ldots, e_l$ and ends at $v_i$. Check whether $u \ast w \in \overline{\mu(q_i)}H_C(Q)$. Since $P$ is
balanced, if the answer is negative for every $i=1,2,\ldots,k$ then there is no such $p$. If the
answer is positive for some $i_0$, then we can take $p = q_{i_0} q$, where $q$ is a loop in $Q$ such
that $u \ast w = \overline{\mu(q_i)}\overline{\mu(q)}$.

\hfill $\square$

\section{Folded ${\cal B}(G)$-graphs}
\label{sec:4}

As we noticed in Section 2, an element $g \in G$ can be represented by its normal form. Ideally, for
a ${\cal B}(G)$-graph $\Gamma$ and its vertex $v$ we would like to reduce the question ``if $g \in
L(\Gamma,v)$'' to the question ``if there exists a loop $p$ in $\Gamma$ at $v$ labeled by $\pi(g)$''.
The idea is that it would be much easier to answer the latter question than the former one, provided
we could construct a completion of $\Gamma$ which recognizes the same language but for which the
above questions are equivalent.

However, certain technical problems arise which makes such a completion graph quite complicated. To
avoid these technical difficulties, we impose a weaker condition on the graph which allows to keep
its structure relatively simple, albeit at the cost of making reading procedure rather involved. To
produce precise statements, we need to introduce some additional terminology first.

\smallskip

Recall that by Theorem \ref{th:1},
$$G_i = \langle G_{i-1}, T_{i-1} \mid t^{-1} C_t t \stackrel{\phi_t}{=} D_t,\ t \in T_{i-1},\ C_t,
D_t \in {\cal C}_{i-1} \rangle.$$
\begin{defn}
\label{de:4.2} Let $\Gamma$ be a finite ${\cal B}(G)$-graph.
Then $\Gamma(i)$ is a subgraph of $\Gamma$ such that if $e \in E(\Gamma(i))$ then either
\begin{enumerate}
\item $\mu(e) = x \in X^\pm$, or
\item $\mu(e) \in T_j,\ j < i$, or
\item $\mu(e) = v,\  v \in C\in \mathcal C_j,\ j < i$.
\end{enumerate}
We call $\Gamma(i)$ the {\em $i$-level graph of $\Gamma$} (for example, the $1$-level graph is the
subgraph of $\Gamma$ which consists only of edges with labels from $X$) and say that $\Gamma$ has
{\em level} $N$ denoted by $l(\Gamma)$ if $N$ is the minimal natural number for which $\Gamma =
\Gamma(N)$.
\end{defn}
Observe that $\Gamma(i)$ may not be connected for some $i < l(\Gamma)$, but all notions introduced
in Section \ref{sec:3} still apply.

As we mentioned in Remark \ref{re:2.2.1}, normal forms are not unique. With that in mind we
introduce a weaker version of standard decomposition for elements of $G$ that still relies on the
normal form theorem for HNN-extensions.
\begin{defn}
To each element $g \in G$ we associate an (infinite) set of words $\Pi(g)$ as follows. If $g \in
G_1$, then $\Pi(g) = \{\pi(g)\}$. Otherwise, assuming
$$\pi(g) = \pi(g_1) w_1(r_1) \cdots \pi(g_k) w_k(r_k) \pi(g_{k+1}) = \pi(g_1) w_1(r_1) \pi(\overline{g}),$$
denote $C = C_{w_1},\ D = D_{w_1}, \phi = \phi_{w_1}$ and put
$$\Pi(g) = \bigcup_{\stackrel{u \in C, v' \in D}{v = \phi(u)}} \Pi(g_1 \ast u) (w_1 v') \Pi(v_{w_1}^{r_1}
\ast v^{-1} \ast v'^{-1} \ast \overline{g}).$$
\end{defn}
This definition does not depend on a particular choice of $\pi(g)$ by Lemma \ref{le:2.2.2}.

\subsection{Weakly folded ${\cal B}(G)$-graphs}

In this subsection we introduce the conditions mentioned above and investigate the properties of
${\cal B}(G)$-graphs that accommodate these conditions.

\begin{defn}
Let $C \in {\cal C}_j$. A $C$-component $Q$ is called \emph{label-maximal}, if for any path $p$ in
$\Gamma$ with $o(p) \in V(Q)$ and $\overline{\mu(p)} \in C$ there exists a path $q$ in $Q$ with
$o(q) = o(p), t(q) \in V(Q)$ and $\overline{\mu(q)} = \overline{\mu(p)}$.
\end{defn}

Let $Q$ be a $C$-component for $C \in {\cal C}_j$. We say that the {\em paths in $Q$ are doubled} if
for a path $q$ in $Q$ labeled by $u \in C$ there exists a path $p$ in $\Gamma$ with $o(p) = o(e),
t(p) = t(e), \mu(p) \in \Pi(u)$.

\begin{defn}
We call a ${\cal B}(G)$-graph $\Gamma$ of level $N + 1,\ N \geqslant 0$ {\em weakly folded} if it
has the following properties:
\begin{itemize}
\item[(P1)] $\Gamma(N)$ is weakly folded,
\item[(P2)] $\Gamma$ is freely folded,
\item[(P3)] each $C$-component of $\Gamma$ is reduced ($C \in {\cal C}_N$),
\item[(P4)] paths in all $C$-components are doubled ($C \in {\cal C}_N$),
\item[(P5)] each $C$-component of $\Gamma$ is label-maximal ($C \in {\cal C}_N$),
\item[(P6)] each $ST_N$-complex of $\Gamma$ is reduced and balanced.
\end{itemize}
In particular, a graph of level $1$ is weakly folded if it is freely folded.
\end{defn}

The following theorem provides statement of conditions and properties mentioned in the beginning
of this section.

\begin{thm}
\label{th:4.1}
The following statements hold for each integer $N \geqslant 1$.
\begin{itemize}
\item[${\bf R}(N)$] Let $\Gamma$ be a weakly folded ${\cal B}(G)$-graph of level $N$. Let $v_1, v_2
\in V(\Gamma)$. Let $g \in G$ be given as a word in the generators. Then there is an algorithm to
check whether there exists a path $p$ in $\Gamma$ with $o(p) = v_1, t(p) = v_2, \mu(p) \in \Pi(g)$.

\item[${\bf E}(N)$] Let $\Gamma$ be a weakly folded ${\cal B}(G)$-graph of level $N$. Let $g \in G$
and let $p$ be a path $\Gamma$ such that $\overline{\mu(p)} = g$. Then there exists a path $q$ in
$\Gamma$ with $o(q) = o(p)$, $t(q) = t(p)$ and $\mu(p) \in \Pi(g)$.

\item[${\bf C}(N)$] Let $\Gamma$ be a ${\cal B}(G)$-graph of level $N$. Let $v_1, v_2 \in V(\Gamma)$.
There is an algorithm to construct a weakly folded ${\cal B}(G)$-graph $\Gamma'$ and a generalized
morphism $\Gamma \to \Gamma'$, $v_1 \mapsto v'_1, v_2 \mapsto v_2'$ with $L(\Gamma, v_1, v_2) =
L(\Gamma', v_1', v_2')$.
\end{itemize}
\end{thm}

We prove Theorem \ref{th:4.1} by induction on $N$. The statements ${\bf R}(1), {\bf E}(1), {\bf C}(1)$
are about Stallings graphs in free groups and, therefore, they are clear. For $N > 1$, the induction
steps for ${\bf R}(N), {\bf E}(N), {\bf C}(N)$ are shown in the following subsections. We will also
establish a slightly stronger version of ${\bf R}(N)$ (coinciding with ${\bf R}(N)$ for $N = 1$).

{\it
\begin{itemize}
\item[${\bf R^+}(N)$] Let $\Gamma$ be a weakly folded ${\cal B}(G)$-graph of level $N$. Let $v_1 \in
V(\Gamma)$. Let $C \in {\cal C}_{N-1}$ and $g \in G_{N}$ be given as a word in the generators. Then
there is an algorithm to check whether there exists an element of $G$ of the form $g u, u \in C$
which is readable from $v_1$ in $\Gamma$. Moreover, there is an algorithm to produce a finite list
of elements $g u_1, g u_2, \ldots$, where $u_i \in C$, and a $C$-component $Q$ in $\Gamma$ such that
each $g u_i$ is readable from $v_1$ in $\Gamma$ along a path terminating at a vertex of $Q$, and,
further, any element of the form $g u, u \in C$ readable from $v_1$ in $\Gamma$ belongs to the coset
$(g u_i) H_C(Q)$.
\end{itemize}}

With a slight abuse of notation, we also call the effective procedures whose existence is provided
in the statements ${\bf R}(N), {\bf R^+}(N), {\bf C}(N)$ by the same letters ${\bf R}(N), {\bf R^+}(N),
{\bf C}(N)$ respectively.

The induction steps will be arranged as follows:
\[\begin{array}{ll}
\cdots & \Rightarrow {\bf E}(N) \Rightarrow {\bf E}(N+1) \Rightarrow \cdots,\\
\cdots & \Rightarrow {\bf R}(N), {\bf E}(N), {\bf C}(N), {\bf R^+}(N) \Rightarrow \\
& \Rightarrow {\bf C}(N+1) \Rightarrow {\bf R}(N+1), {\bf R^+}(N+1) \Rightarrow \cdots.
\end{array}
\]

We would like to note that dependence of the reading procedure ${\bf R}(N+1)$ on the folding procedure
${\bf C}(N+1)$ at the same level is not strictly necessary, as ${\bf R}(N+1)$ can be also deduced from
${\bf R^+}(N)$. However, neither way provides an advantage over the other one in terms of
computational complexity, so for the sake of making shorter argument, we choose the former way.

\subsection{Induction step for ${\bf E}(\cdot)$}

Suppose that the statement ${\bf E}(N)$ holds. In this subsection we prove ${\bf E}(N+1)$.

Let $q$ be a path in a weakly folded graph $\Gamma$ of level $N + 1$, and $\overline{\mu(q)} = g$.
Our goal is to show that there is also a path $p$ from $o(q)$ to $t(q)$ that reads an element from
$\Pi(g)$.

Let
$$\pi(g) = \pi(g_1) w_1(r_1) \cdots \pi(g_k) w_k(r_k) \pi(g_{k+1}).$$
At first we show that we can assume
$$\mu(q) = h_1 (w_1 v_1) h_2 \cdots h_k (w_k v_k) h_{k+1},$$
where $v_i \in D_{w_i}$ for $i \in [1, k]$. Indeed, by (P4) we can assume without loss of generality
that
$$q = d_1 e_1 d_2 \cdots d_l e_l d_{l+1},$$
where all $d_i$ are paths in $\Gamma(N)$ and $e_i$ are maximal (with respect to inclusion)
$ST_N$-paths of maximal height. Denoting $t_i = \mu(e_i)$ and $\mu(d_i) = h'_i, \overline{h'_i} \in
G_N$ we have
$$\mu(q) = h_1' t_1 h_2' t_2 \cdots h'_l t_l h'_{l+1}.$$
If for every $i \in [1, l-1]$ we have $D_{t_i} \neq C_{t_{i+1}}$, or $d_{i+1} \notin D_{t_i}$, then
by Lemma \ref{le:2.2.2}, $k = l$ and $w_i t_i^{-1} \in C_{w_i} = C_{t_i}$ for $i \in [1, l]$.
Suppose $l \neq k$, therefore, there exists $i$ such that $D_{t_i} = C_{t_{i+1}}$ and $d_{i+1} \in
D_{t_i}$. Then by (P5) we can assume that $d_{i+1}$ is a path inside of the $D_{t_i}$-component of
$t(t_i)$, and $t_i d_{i+1} t_{i+1}$ can be replaced either by a path in $C_{t_i}$-component of
$o(t_i)$ (if $C_{t_i} = D_{t_{i+1}}$), or by an $ST_N$-path (otherwise). After such a replacement,
the number $l$ is reduced either by $2$, or by $1$, respectively.

Repeating this argument, we ultimately force $l = k$ and $w_i t_i^{-1} \in C_{w_i} = C_{t_i}$ for
$i \in [1, k]$.

We proceed by induction on $k$. If $k = 0$, by (P4) we assume that $q$ is a path in $\Gamma(N)$ and
by (P1) and ${\bf E}(N)$ there exists a path $p$ in $\Gamma(N)$ such that $\mu(p) \in \Pi(g)$.

If $k > 0$ then we assume that we have established the statement for elements with less than $k$
syllables and denote the piece of $q$ corresponding to $h_1$ by $q_1$. By (P4), we assume $q_1$ to
be a path in $\Gamma(N)$, which, together with (P1), enables us to use ${\bf E}(N)$. Therefore, we
may assume that $\mu(q_1) \in \Pi(h_1)$. Now reading $\Pi(g)$ is reduced to reading $\Pi(h_2 w_2
\cdots w_k h_{k+1})$, which has at most $k - 1$ syllables.

\begin{rmk}
\label{rem:subdivision}
Let $\pi(g)$ be as described in the beginning of the subsection. Suppose some path $p$ reads $\pi'(g)
\in \Pi(g)$. Then, given $l \in [1, k]$, we cannot guarantee that $g_1 \circ w_1(r_1) \circ \cdots
\circ g_l$ is readable along an initial subpath of $p$. Nevertheless, a weaker but still useful
statement can be made: $p$ can be subdivided into subpaths $p = p_1 p_2 p_3$ so that
\[\begin{array}{rcl}
\overline{\mu(p_1)} & = & g_1 \ast w_1(r_1) \ast \cdots \ast g_l \ast u_l,\\
\overline{\mu(p_2)} & = & u'_l \ast w_l \ast v'_l,\\
\overline{\mu(p_3)} & = & v_l \ast g_{l+1} \ast \cdots \ast w_k(r_k) \ast g_{k+1},\\
\end{array}
\]
where $u_l, u'_l \in C_{w_l}, v_l, v'_l \in D_{w_l}$ and $p_2$ is an $ST_N$-path.
\end{rmk}

\subsection{Induction step for ${\bf C}(\cdot)$}

Now we present a procedure that, given a ${\cal B}(G)$-graph $\Gamma$, allows to obtain a weakly
folded graph that recognizes the same language. We proceed by induction on the level of $\Gamma$,
that is, we assume ${\bf R}(N), {\bf R^+}(N), {\bf E}(N), {\bf C}(N)$ to hold, and prove
${\bf C}(N+1)$.

Let $\Gamma$ be a graph of level $N + 1$. By ${\bf C}(N)$ we can assume that $\Gamma(N)$ is weakly
folded, thus gaining (P1). Note also that for graphs of level greater than $1$ (P2) follows from (P1),
(P3) and (P6).

Throughout this section, when we write $L(\Gamma) = L(\Gamma')$ we imply existence of a generalized
morphism $\Gamma \to \Gamma'$ and we mean the equality $L(\Gamma, v_1, v_2) = L(\Gamma', v'_1, v'_2)$
for all $v_1, v_2 \in V(\Gamma)$ and $v_1 \mapsto v'_1,\ v_2 \mapsto v_2'$.

\begin{lem}
\label{le:4.20}
Let $\Gamma$ be a finite ${\cal B}(G)$-graph of level $N + 1$ with the property (P1). Then using
free foldings, simple reductions and the algorithm ${\bf C}(N)$ finitely many times, we can obtain a
graph $\Gamma'$ with $L(\Gamma) = L(\Gamma')$ for which (P1) and (P3) hold.
\end{lem}
{\it Proof.} Using finitely many simple reductions we can produce a graph $\Gamma_1$, where all
$C$-components of $\Gamma$ ($C \in {\cal C}_N$) are reduced. Note that $\Gamma_1(N)$ may not be weakly
folded anymore. This can only happen if some vertices in a $C$-component were identified, that is, if
$$\sum_Q |V(Q)| > \sum_{Q_1} |V(Q_1)|,$$
where $Q_1$ is a $C$-component of $\Gamma_1$ obtained from $Q$. By the induction hypothesis we can
transform $\Gamma_1$ into a graph $\Gamma_2$ such that $\Gamma_2(N)$ is weakly folded. Let components
$Q_1$ transform into $Q_2$. Since edges labeled by $u \in C$ of non-maximal height belong to weakly
folded $\Gamma_2(N)$, the only way $Q_2$ may cease to be reduced is when
$\mathrm{proj}_{u_\mathrm{max}} H_C(Q_1)$ gets changed, where $u_\mathrm{max}$ is the generator of
$C$ of maximal height and $\mathrm{proj}_{u_\mathrm{max}}$ stands for the projection onto the cyclic
group $\langle u_\mathrm{max} \rangle$. So, one of components $Q_2$ may be not reduced only if two
vertices in (coinciding or distinct) $C$-components $Q_1$ have been identified while transforming
$\Gamma_1$ into $\Gamma_2$, which increases $H_C(Q_1)$. Since finitely generated free abelian groups
are noetherian, this can only happen finitely many times.

Therefore, in a finite number of steps we obtain a graph $\Gamma'$ for which (P1) and (P3) hold.
$L(\Gamma) = L(\Gamma')$ by Lemma \ref{le:3.2.2} and Lemma \ref{le:3.4.1}.

\hfill $\square$

\begin{rmk}
Note that strictly speaking, the algorithmic step in Lemma \ref{le:4.20} is redundant since (P3)
follows from (P1) and (P4). However, in certain cases obtaining (P3) prior to (P4) provides lower
computational complexity.
\end{rmk}

\begin{lem}
\label{le:4.21}
Let $\Gamma$ be a finite ${\cal B}(G)$-graph of level $N + 1$ with the properties (P1) and (P3).
Then adding finitely many vertices and edges, and using free foldings, simple reductions and the
algorithm ${\bf C}(N)$ finitely many times, we can obtain a graph $\Gamma'$ with $L(\Gamma) =
L(\Gamma')$ for which (P1), (P3) and (P4) hold.
\end{lem}
{\it Proof} Let $Q$ be a $C$-component for $C \in {\cal C}_N$. By ${\bf R}(N)$, for each edge $e
\in Q$ check whether an element of $\Pi(\mu(e))$ is readable from $o(e)$ to $t(e)$ in $\Gamma(N)$.
In the case it is not readable (which, by (P1), implies that $\overline{\mu(e)}$ is of maximal
height in $C$), add a path $q$ from $o(e)$ to $t(e)$ with $\mu(q) \in \Pi(\mu(e))$. Note that this
does not change the language defined by the graph $\Gamma$.

Thus, having added finitely many edges and vertices we provide that all paths are doubled in
$\Gamma(N)$. This may have broken (P1). Using Lemma \ref{le:4.20}, obtain (P1) and (P3). Since the
procedure in Lemma \ref{le:4.20} does not produce new edges labeled by elements of maximal height
in $C$, the property (P4) still holds.

\hfill $\square$

In order to obtain (P5), we distinguish two cases: when the rank of $C$ is equal to $1$ and when
the rank of $C$ is bigger than $1$, where $C$ is one of the conjugated centralizers $C \in {\cal
C}_N$. The following lemma is a technical statement we use to treat the former case.

\begin{lem}
\label{4.loop}
Let $\Delta$ be a weakly folded ${\cal B}(G)$-graph of level $N,\ v \in V(\Delta)$. Let $g \in G_N$
and paths $p_i,\ i = 1, 2, \ldots$ be such that $o(p_i) = v,\ \overline{\mu(p_i)} = v_i \ast g \ast
u_i$, where $v_i \in C \in {\cal C}_{\leqslant N-1},\ u_i \in D \in {\cal C}_{\leqslant N-1}$. Then
\begin{enumerate}
\item[(a)] if $[g^{-1} C g, D] \neq 1$ the there exist $i_0, i_1$ such that $v_{i_0} \ast v_{i_1}^{-1}
\in H_C(v)$,
\item[(b)] if $[g^{-1} C g, D] = 1$ then there exists a subsequence $\{i_j\},\ j = 1, 2, \ldots$ of
$\{i\}$ such that all $p_{i_j}$ terminate at the same vertex $v_t$, and there exist loops $l_j$ in
the $D'$-component of $v_t$ such that $\overline{\mu(p_{i_j} p_{i_{j+1}}^{-1} l_j)} = 1$, where $D'$
is an abelian subgroup of $G$ of maximal height such that $D \subseteq D' \in {\cal C}_{\leqslant
N-1}$.
\end{enumerate}
\end{lem}
{\it Proof.} (b) is obvious. Indeed, denoting $g^{-1} \ast v_i \ast g = u'_i$, we have
$\overline{\mu(p_i)} = g \ast u'_i \ast u_i$, where $u_i, u'_i \in D'$. Then the subsequence $i_j$,
such that $p_{i_j}$ terminate at the same vertex, delivers the statement.

\smallskip

We prove (a) by reducing it to an analogous statement about $G_{N-1}$ and shifts $v'_i, u'_i$ that
are elements of centralizers that belong to ${\cal C}_{< N-1}$.

Denote the generator of $C$ of maximal height by $v_{\mathrm{max}}$ and put $v_i =
v_{\mathrm{max}}^{k_i} \ast v_i'$, where $v_i' \in C' \leqslant C,\ C' \in {\cal C}_{< N-1}$. In
particular, $v_i = v_i'$ if $C \in {\cal C}_{< N-1}$. Similarly, denote the generator of $D$ of
maximal height by $u_{\mathrm{max}}$ and put $u_i = u_{\mathrm{max}}^{m_i} \ast u_i'$, where $u_i'
\in D' \leqslant D,\ D' \in {\cal C}_{< N-1}$. If the set $\{k_i \mid i = 1, 2, \ldots \}$ is
bounded, passing to a subsequence we may assume that $v_i = v_{\mathrm{max}}^k \ast v_i'$. Similarly,
if the set $\{m_i \mid i = 1, 2, \ldots \}$ is bounded, we may assume that $u_i =
u_{\mathrm{max}}^m \ast u_i'$. Using Lemma \ref{le:2.2.1}, and, if necessary, passing to a
subsequence, we conclude that all elements $v_i \ast g \ast u_i$ have one of the following four forms:

\[\mathrm{{\bf Case\ 1}\ } \left\{\begin{array}{l}
v_i \ast g \ast u_i = v'_i \ast [v_{\mathrm{max}}^k \ast g \ast u_{\mathrm{max}}^m] \ast u'_i,\
\mathrm{or}\\
v_i \ast g \ast u_i = v'_i \ast [v_{\mathrm{max}}^k \ast g \ast u_{\mathrm{max}}^M] \circ
u_{\mathrm{max}}^{m_i - M} \ast u'_i,
\end{array}\right.
\]
$$\mathrm{or}$$
\[\mathrm{{\bf Case\ 2}\ } \left\{\begin{array}{l}
v_i \ast g \ast u_i = (v'_i \ast v_{\mathrm{max}}^{k_i-M}) \circ [v_{\mathrm{max}}^M \ast g \ast
u_{\mathrm{max}}^m] \ast u'_i,\ \mathrm{or}\\
v_i \ast g \ast u_i = (v'_i \ast v_{\mathrm{max}}^{k_i-M}) \circ [v_{\mathrm{max}}^M \ast g \ast
u_{\mathrm{max}}^M] \circ (u_{\mathrm{max}}^{m_i-M} \ast u'_i).
\end{array}\right.
\]

Consider the cases above.

\smallskip

{\bf Case 1.} For simplicity, assume the bracketed part to be equal to $g$. Let
$$\pi(g) = \pi(g_1)\ w_1(r_1)\ \pi(g_2)\ w_2(r_2)\ \cdots\  w_k(r_k)\  \pi(g_{k+1}),$$
where $w_i \in ST_{N-1}$.

Suppose $[g_1^{-1} C g_1, C_{w_1}] = 1$ and $k > 0$. Then, since $v'_i \in G_{N-1}$, by ${\bf E}(N)$
(see Remark \ref{rem:subdivision}), the paths $p_i$ can be subdivided into subpaths $q'_i\ q_i = p_i$
that read
$$\overline{\mu(q'_i)} = v'_i \ast g_1 \ast w_1 \ast u''_i,\ \ \overline{\mu(q_i)} = {u''_i}^{-1}
\ast g_2 \ast w_{2} \ast \cdots \ast g_{k+1} \ast u_i,$$
where $u''_i \in D_{w_1}$. Then, by (b) applied to $q_i'$, we pass to a subsequence of $q_i$'s and
the word $g' = g_2 \circ w_2(r_2) \circ \cdots \circ g_{k+1}$ of smaller level $N$ syllable length
than $g$. Either $g'$ has level $N$ syllable length $0$, or greater than $0$ (that is, $g \in G_N -
G_{N-1}$), and $u''_i \in G_{N-1}$.

Hence, repeating the argument, we eventually may assume that $k > 0$ and $[g_1^{-1} C g_1, C_{w_1}]
\neq 1$, or $k = 0$ and, by the assumption of (a), $[g_1^{-1} C g_1, D] \neq 1$. By ${\bf E}(N)$ and
Remark \ref{rem:subdivision}, there exist a collection of paths $q_i$ in $\Delta(N-1)$ reading
$\overline{\mu(q_i)} = v_i' \ast g_1 \ast u_i''$, where either $u_i'' \in C_{w_1}$ if $k > 0$, or
$u_i'' \in D'' \in {\cal C}_{< N-1},\ [D'', D] = 1$, if $k = 0$. The statement is thus reduced to
graphs of lower level. Note that if the rank of $C$ is equal to $1$ then this case delivers the
statement of the lemma, since in this event $v'_i = 1$.

\smallskip

{\bf Case 2.} Show further that in this case we still may assume (concatenating $p_i$ with a loop
at $v$) that the set $\{k_i - M \mid i = 1, 2, \ldots \}$ is bounded. More precisely, show that the
$C$-component $Q$ of $v$ contains a loop labeled by $v_{\mathrm{max}}^{k_0} \ast v'$, where $v' \in
C$ is of non-maximal height.

Indeed, $\pi(v_i \ast g \ast u_i)$ begins with $\pi(v'_i \ast v_{\mathrm{max}}^{k_i - M})$, therefore,
the loop in question exists by ${\bf E}(N)$ and ${\bf E}(N-1)$ since $\Delta$ contains only a finite
number of $C'$-components.

Denote loop at $v$ reading $v_{\mathrm{max}}^{k_0} \ast v'$ by $q$. Then with an appropriate choice
of powers $l_i$ we can make $q^{l_i} p_i$ satisfy {\bf Case 1}.

\hfill $\square$

\begin{lem}
\label{le:4.22}
Let $\Delta$ be a weakly folded ${\cal B}(G)$-graph of level $N$ and let $v \in V(\Delta)$. Let $u
\in G_N$ be such that $\pi(u^2) = \pi(u)^2$. Then there is an algorithm to construct the label-maximal
$\langle u \rangle$-component at $v$.
\end{lem}
{\it Proof.} In other words, we need to find all end vertices of paths $p$ that read $u^k$ starting
from $v$.

Let $\pi(u) = \pi(h_1)\ w_1(r_1)\ \pi(h_2)\ w_2(r_2)\ \cdots\ w_k(r_k)\ \pi(h_{k+1})$. By
${\bf R^+}(N)$, we consecutively attempt to read $\Pi(u),\ \Pi(u)\ h_1\ v',\ \Pi(u^2),\ \Pi(u^2 h_1)\
v', \ldots$, where $v' \in C = C_{w_1}$, starting from $v$ in $\Delta$. We stop if either
\begin{enumerate}
\item[(a)] some $\Pi(u^N h_1) v'$ is not readable for any $v' \in C$, or
\item[(b)] we discover that some $\Pi(u^N)$ is readable as a label of a loop at $v$.
\end{enumerate}

If (a) occurs, by ${\bf E}(N)$ and Remark \ref{rem:subdivision} no loop labeled by a power of $u$
is readable at the vertex $v$, so it is enough to check which of $\Pi(u^{\pm 1}), \ldots,
\Pi(u^{\pm N})$ are readable as labels of paths from $v$.

If (b) occurs, then for any $u^M$ readable starting from $v,\ u^{M \pm N}$ is also readable with
the same terminal vertex, so it is enough to check which of $\Pi(u^{\pm 1}), \ldots, \Pi(u^{\pm N})$
are readable as labels of paths from $v$.

Now we show that it takes only finitely many steps for either one of (a) and (b) to occur. Either
at some point $\Pi(u^N h_1) v'$ is readable for no $v' \in C$, or endpoints of the
corresponding paths meet some $C$-component $Q$ arbitrarily many times.

The former case is exactly (a). In the latter case, show that there exists a loop at $v$ reading a
power of $u$.

Note that if paths $p_i$ ($i = 1, 2, \ldots$) are such that $o(p_i) = v,\ t(p_i) \in Q,\
\overline{\mu(p_i)} = u^{k_i} \ast h_1 \ast v_i$, then concatenating with a path in $Q$ we may
assume that $t(p_i) = t(p_j)$ for all $i, j$. Then by ${\bf E}(N)$, each $p_i$ contains a terminal
subpath $q_i$ with $\overline{\mu(q_i)} = u_i \ast h_{k+1} \ast h_1 \ast v_i$, where $u_i \in
D_{w_k}$. Since $\pi(u^2) = \pi(u)^2$, the paths $q_i^{-1}$ satisfy conditions of Lemma \ref{4.loop}
(a), therefore, for some loop $q$ in $Q$ we have $\overline{\mu(p_{i_0} q p^{-1}_{i_1})} = u^l$,
where $l = k_{i_0} - k_{i_1}$.

\hfill $\square$

\begin{rmk}
The technical difficulty in solving power problem for $u$ in $L(\Delta)$ by direct inspection of
$\Delta$ is caused by the fact that there is no clear connection between $\Pi(u)$ and $\Pi(u^2)$,
which makes it theoretically possible to have a weakly folded graph $\Delta$, in which $\Pi(u^M)$
can be read as a label of a loop at a vertex $v$, but $\Pi(u^{< M})$ cannot be read as a label of
a path from $v$. We would like to point out that this is essentially the same difficulty that one
faces trying to solve the membership problem for $\mathbb Z^n$-free groups using the machinery of
\cite{KMW}, as we mentioned in Section \ref{sub:structure} on page \pageref{difficuly}.
\end{rmk}

\begin{lem}
\label{le:4.23}
Let $\Delta$ be a weakly folded ${\cal B}(G)$-graph of level $N$ and let $v \in V(\Delta)$. Let $C
\in {\cal C}_N$ be such that the rank of $C$ is greater than $1$. Then there is an algorithm to
construct the label maximal $C$-component at $v$.
\end{lem}
{\it Proof.} In other words, we need to find all end vertices of paths whose labels are equal to
$u = u_{{\mathrm{max}}}^k \ast u'$, where $u_{\mathrm{max}}$ is a generator of $C$ of maximal height
and $u'$ is an element of less than maximal height. To do that it is enough to inspect the
$ST_{N-1}$-complex of $v$ using Lemma \ref{le:3.5.4}.

\hfill $\square$

\begin{lem}
\label{le:4.24}
Let $\Gamma$ be a finite ${\cal B}(G)$-graph of level $N + 1$ with the properties (P1), (P3) and (P4).
Then adding finitely many vertices and edges, using free foldings, simple reductions and the algorithm
${\bf C}(N)$ finitely many times, we can obtain a graph $\Gamma'$ with $L(\Gamma) = L(\Gamma')$ for
which (P1), (P3)--(P5) hold.
\end{lem}
{\it Proof.} To acquire (P5), we apply either Lemma \ref{le:4.22}, or Lemma \ref{le:4.23} to all $C
\in {\cal C}_N$ and all $v \in V(\Gamma(N))$. Note that since $\Gamma(N)$ is weakly folded this
procedure does not break (P1), (P3) and (P4).

\hfill $\square$

\begin{lem}
\label{le:4.25}
Let $\Gamma$ be a finite ${\cal B}(G)$-graph of level $N + 1$ with the properties (P1), (P3)--(P5).
Then adding finitely many vertices and edges, using free foldings, simple reductions, loop translations,
$ST_N$-foldings and the algorithm ${\bf C}(N)$ finitely many times we can obtain a graph $\Gamma'$ with
$L(\Gamma) = L(\Gamma')$ for which (P1)--(P6) hold.
\end{lem}
{\it Proof.} Let $P$ be an $ST_N$-complex. As shown before, using finitely many loop translations
and $ST_N$-foldings, we can transform $\Gamma$ into a graph $\Gamma_1$ such that the $ST_N$-complex
$P_1$ corresponding to $P$ is a reduced balanced complex. Then, using previous lemmas, we transform
$\Gamma_1$ into a graph $\Gamma_2$ for which (P1), (P3)--(P5) hold. The $ST_N$-complex $P_2$
corresponding to $P_1$ may cease to be reduced in the only case: some $C$-components of $\Gamma_1$
($C \in {\cal C}_N$) got glued together. Since the number of $C$-components is finite, repeating
this procedure we may assume $P_2$ to be reduced. Then $P_2$ may cease to be balanced only if some
$C$-components got changed in the process, increasing $H_C(Q)$. This can only happen a finite number
of times, since $H_C(Q)$ is a subgroup of a noetherian group.

Eventually, we obtain a graph $\Gamma'$ for which (P1), (P3)--(P6) (and therefore (P2)) hold.
$L(\Gamma) = L(\Gamma')$ by Lemmas \ref{le:3.5.1},\ \ref{le:3.5.2} and \ref{le:4.20}--\ref{le:4.24}.

\hfill $\square$

\begin{thm}
Let $\Gamma$ be a finite ${\cal B}(G)$-graph of level $N + 1$. Then adding finitely many vertices
and edges, using free foldings, simple reductions, loop translations, $ST_N$-foldings and the
algorithm ${\bf C}(N)$ finitely many times we can obtain a weakly folded graph $\Gamma'$ with
$L(\Gamma) = L(\Gamma')$.
\end{thm}
{\it Proof.} Immediate from Lemmas \ref{le:4.20}--\ref{le:4.25}.

\hfill $\square$

Note that all the operations listed in the theorem above are either graph morphisms, or generalized
graph morphisms, so for each vertex $v$ of $\Gamma$ there is a corresponding vertex $v'$ in $\Gamma'$.

\subsection{Induction step for ${\bf R}(\cdot)$ and ${\bf R^+}(\cdot)$}

Assume that the statements ${\bf R}(N),\ {\bf E}(N),\ {\bf C}(N+1)$ hold. In this subsection we prove
${\bf R}(N+1)$ and ${\bf R^+}(N+1)$.

Fix a weakly folded graph $\Gamma$ of level $N + 1$. Suppose $g \in G_{N+1}$. Build a graph $\Gamma'$
of level $N + 1$ by attaching a path $p_g$ reading $g$ to $\Gamma$ at a vertex $v_1$, that is, such
that $o(p_g) = v_1,\ \mu(p_g) \in \Pi(g)$. Next, by ${\bf C}(N+1)$ construct a weakly folded graph
$\Gamma''$ with $L(\Gamma'') = L(\Gamma')$. Denote by $v''$ the vertex in $\Gamma''$ that corresponds
to $t(p_{g_1})$ in $\Gamma'$. The question whether $\Pi(g)$ is readable from $v_1$ to $v_2$ in
$\Gamma$ is now reduced to the question whether $v''$ corresponds to the vertex $v_2$ of $\Gamma$.

To establish ${\bf R^+}(N+1)$ it is enough to inspect the $C$-component ($C \in {\cal C}_N$) of $v''$.
More exactly, we note that $\Gamma''$ defines the same subgroup of $G$ as $\Gamma$, so no vertices
of (distinct or coinciding) $C$-components of $\Gamma$ could be identified while obtaining $\Gamma''$.
Suppose a path $g u,\ u \in C$ is readable in $\Gamma$ from $v_1$. Then it is also readable in
$\Gamma''$ along the path $p_{g u}$ from the image of $v_1$. By (P4), $v''$ and $t(p_{g u})$ belong
to the same $C$-component $Q''$ of $\Gamma''$. Denote by $Q$ the $C$-component of $\Gamma$ that
corresponds to $Q''$ in $\Gamma''$ and let $V(Q)$ be the finite set of vertices of $Q$. Denote by
$V''(Q'')$ the set of vertices of $Q''$ that correspond to $V(Q)$. Denote by $u_1, u_2, \ldots$
labels of arbitrary paths inside $Q''$ from $v''$ to vertices in $V''(Q'')$. Then the elements $g \ast
u_1, g \ast u_2, \ldots$ and the $C$-component $Q$ deliver the statement ${\bf R^+}(N+1)$.

\section{Results and conclusions}
\label{sec:5}

Before giving some corollaries of Theorem \ref{th:4.1} we discuss applicability of the obtained
results.

\subsection{Note on effectiveness}
\label{ssec:5.1}

The main question that arises is the following. Suppose the input is a presentation of group a $G$ and
the fact that it acts freely and regularly on {\em some} $\mathbb Z^n$-tree. Can we solve the uniform
subgroup membership problem in $G$ using the techniques developed in the present work? Strictly
speaking, the answer is unfortunately negative (at least, to the date). Indeed, Theorem \ref{th:4.1}
gives an effective construction of a weakly folded ${\cal B}(G)$-graph and subsequently means to
read words as labels of paths in the graph, thus solving the uniform membership problem for a fixed
ambient $G$ (see Theorem \ref{th:5.1} below). More exactly, it proves existence of a decision
algorithm but only gives a way to obtain such an algorithm if the data in Theorem \ref{th:4.1} is
given as a part of the input. However, if we are concerned with uniformity over different groups $G$,
to build and use weakly folded graphs, we have to address certain issues.

To elaborate, suppose first that the chain of HNN-extensions for $G$ as described in Theorem
\ref{th:main2} is given as a part of the input (we say in this case that $G$ is {\em given by an
effective HNN-chain}). Can one then assume that the Lyndon length function of any element of $G$ can
be computed effectively? The answer is positive. Indeed, as one can see from the proof of Theorem
\ref{th:main2} and Theorem \ref{main4} (see \cite{KMRS}), an effective HNN-chain provides an
algorithm to embed $G$ into a partial group of infinite words, and, moreover, the images of generators
are recursive as functions defined on segments of $\mathbb Z^n$. Using standard diagonal argument,
one can effectively find $com(f,g)$ if $f, g$ are given as recursive $\mathbb Z^n$-words. Then it is
easy to see that the embedding, the group operation and the Lyndon length function are effective.

Therefore, the question is reduced to whether one can algorithmically obtain an effective HNN-chain
given a group presentation and information that the group acts freely and regularly on a
$\mathbb Z^n$-tree. To the date, such an algorithm is not known. However, a number of researchers
express belief that a Makanin--Razborov type machinery can be employed to give effective versions
of Theorems \ref{chis0}, \ref{chis} and ultimately Theorem \ref{th:main2}.

Note that similar questions arise even regarding solution of the Word Problem in $G$. Indeed, a
simple solution that follows from Bass--Serre theory (presented, for example, in \cite{MarR}),
relies on constructing a graph of groups based on the group action on a $\mathbb Z^n$-tree. While this
proves existence of a decision algorithm for a fixed $G$, it does not suggest any means to produce
such an algorithm given a presentation of $G$ and information that it acts freely on some
$\mathbb Z^n$-tree.

Finally, we would like to point out that requiring effectiveness of the Lyndon length function (or
of the embedding into infinite words) does not seem unreasonable. Indeed, it is hard to expect that
information of existence of an action on a $\mathbb Z^n$-tree comes ad hoc, without robust relation
to the presentation of $G$. This thought is reinforced by the fact that the most powerful to the
date tool to produce examples of $\mathbb Z^n$-free groups, Theorem \ref{main4}, implies
effective embedding into infinite words, as mentioned above.

\subsection{Solutions to the Subgroup Membership and Power Problems}
\label{ssec:5.2}

The following theorem states that the uniform membership problem is decidable in $\mathbb Z^n$-free
groups.
\begin{thm}
\label{th:5.1}
Let $G$ be a finitely generated $\mathbb Z^n$-free group. There exists an algorithm that, given $g,
h_1, h_2, \ldots, h_k \in G$ presented as words in the generators of $G$, decides whether $g \in
\langle h_1, h_2, \ldots, h_k \rangle$.
\end{thm}
{\it Proof.} Represent $G$ as a chain of HNN-extensions as described in Theorem \ref{th:1}. Find
normal forms $\pi(g), \pi(h_1), \pi(h_2), \ldots, \pi(h_k)$. Consider the graph $\Gamma_0$ which is
a bouquet of $k$ loops $q_1, \ldots, q_k$ labeled by $\mu(q_i) = \pi(h_i),\ i \in [1,k]$. Apply
${\bf C}(N)$ (Theorem \ref{th:4.1}) for $N = K(G)$ to build a weakly folded version of $\Gamma_0$,
denoted $\Gamma$. Using ${\bf R}(N)$ (Theorem \ref{th:4.1}) decide whether an element of $\Pi(g)$
is readable as a label of a loop in $\Gamma$. By ${\bf E}(N)$ (Theorem \ref{th:4.1}), the positive
answer implies $g \in \langle h_1, h_2, \ldots, h_k \rangle$, while the negative answer implies
$g \notin \langle h_1, h_2, \ldots, h_k \rangle$.

\hfill $\square$

The following theorem states that the power problem is decidable in $\mathbb Z^n$-free groups.

\begin{thm}
\label{th:5.2}
Let $G$ be a finitely generated $\mathbb Z^n$-free group. There exists an algorithm that, given $g,
h_1, h_2, \ldots, h_k \in G$ presented as words in generators of $G$, decides whether there exists
an integer $l \neq 0$ such that $g^l \in \langle h_1, h_2, \ldots, h_k \rangle$.
\end{thm}
{\it Proof.} Represent $G$ as a chain of HNN-extensions as described in Theorem \ref{th:1}. Find $c
= com(g, g^{-1}) \in G$, so that $g_0 = c^{-1} \ast g \ast c$ is cyclically reduced. Then by Lemma
\ref{le:2.2.5} we can find $d \in G$ such that $g_1 = d^{-1} \ast g_0 \ast d = d^{-1} \ast c^{-1}
\ast g \ast c \ast d$ is either cyclically normal, or $g_1 = u \ast w,\ u \in C_w = D_w$.

Consider the graph $\Gamma_0$ that is a bouquet of $k$ loops $q_1, \ldots, q_k$ labeled by $\mu(q_i)
= \pi(h_i),\ i \in [1,k]$ and a path $q$ attached to origin of the bouquet such that $\mu(q) =
\pi(c \ast d)$. By ${\bf C}(N)$ (Theorem \ref{th:4.1}) build a weakly folded version of $\Gamma_0$,
denoted $\Gamma$. Then using either Lemma \ref{le:4.22} (in the former case), or Lemma \ref{le:4.23}
(in the latter case) check whether a loop labeled by a power of $g_1$ is readable from $t(q)$ in
$\Gamma$. The positive answer implies that
$$g_1^l \in \langle d^{-1} c^{-1} h_1 c d,\ d^{-1} c^{-1} h_2 c d,\ \ldots, d^{-1} c^{-1} h_k c d
\rangle,$$
which equivalent to $g^l \in \langle h_1, h_2, \ldots, h_k \rangle$. The negative answer implies that
there is no $l$ such that $g^l \in \langle h_1, h_2, \ldots, h_k \rangle$.

\hfill $\square$

As discussed in Section \ref{ssec:5.1}, Theorem \ref{th:5.1} and Theorem \ref{th:5.2} state existence
of decision algorithms but do not provide means to find such algorithms. The following theorem
provides a way to make the group $G$ a part of the input, although at the cost of making an extra
assumption on $G$.

\begin{thm}
\label{th:5.3}
There exists an algorithm that given
\begin{enumerate}
\item a group $G$ as an effective HNN-chain,
\item elements $g, h_1, h_2, \ldots, h_k \in G$ presented as words in the generators of $G$,
\end{enumerate}
decides whether $g \in \langle h_1, h_2, \ldots, h_k \rangle$ and whether there exists an integer
$l \neq 0$ such that $g^l \in \langle h_1, h_2, \ldots, h_k\rangle$.
\end{thm}
{\it Proof.} Follows from the discussion in Section \ref{ssec:5.1} and Theorems \ref{th:5.1} and
\ref{th:5.2}.

\hfill $\square$

\bibliographystyle{model1-num-names}
\bibliography{Znbib}

\end{document}